\documentclass
{siamltex}
\usepackage{amsfonts}
\usepackage{amsmath}
\usepackage{amssymb}
\usepackage{array}
\usepackage{stmaryrd}
\usepackage{graphicx}
\usepackage{color}
\usepackage{tikz}
\usepackage{enumitem}
\usepackage{bbm}

\usepackage[T3,T1]{fontenc}
\DeclareSymbolFont{tipa}{T3}{cmr}{m}{n}
\DeclareMathAccent{\invbreve}{\mathalpha}{tipa}{16}

\newtheorem{remark}[theorem]{\it Remark}
\newtheorem{theorema}{Theorem}
\newtheorem{theoremaa}{Theorem}

\numberwithin{equation}{section}

\newcommand{\pt}{\partial}

\newcommand {\beq} {\begin{equation}}
\newcommand {\eeq} {\end{equation}}
\renewcommand{\sim}{\simeq}

\newcommand {\B} {{\mathcal B}}

\newcommand {\U} {{\mathcal U}}
\newcommand {\W} {{\mathcal W}}
\newcommand {\N} {{\mathcal N}}

\newcommand{\LL}{{\mathcal L}}
\renewcommand{\AA}{{\mathcal A}}

\newcommand{\RR}{{\mathcal R}}
\newcommand{\R}{\mathbb{R}}

\newcommand{\CCC}{{\lambda}}

\definecolor{pass}{rgb}{0,0,1}
\definecolor{red}{rgb}{0,0,0.7}
\definecolor{blue}{rgb}{0,0,0.7}
\definecolor{green}{rgb}{0,0.6,0}
\begin{document}

\title{Error analysis for time-fractional\\ semilinear parabolic equations
using\\ upper and lower solutions%
\thanks{The author acknowledges support from Science Foundation Ireland Grant SFI/12/IA/1683.}%
}

\author{Natalia Kopteva%
\thanks{Department of Mathematics and Statistics, University of Limerick, Limerick, Ireland ({\tt natalia.kopteva@ul.ie}).}%
 }

\date{}
\maketitle
\begin{abstract}
A semilinear initial-boundary value problem with a Caputo time derivative of fractional order $\alpha\in(0,1)$ is considered,
 solutions of which typically exhibit a singular behaviour at an initial time.
For L1-type discretizations of this problem, we
employ the method of upper and lower solutions
to obtain sharp pointwise-in-time error bounds
on quasi-graded temporal meshes with arbitrary degree of grading.
In particular, those results imply that milder (compared to the optimal) grading yields the optimal convergence rate $2-\alpha$ in positive time, while quasi-uniform temporal meshes yield first-order convergence in positive time.
Furthermore, under appropriate conditions on the nonlinearity,
the method of upper and lower solutions immediately implies  that, similarly to the exact solutions, the computed solutions lie within a certain range.
Semi-discretizations in time and full discretizations using finite differences and finite elements in space are addressed.
The theoretical findings are illustrated by numerical experiments.
\end{abstract}

\section{Introduction}

The method of upper and lower solutions is a very elegant technique frequently used in the analysis of semilinear parabolic and elliptic equations \cite{amann,evans,pao_book},
as well as their discretizations 
\cite{pao_numerical,pao_numerical1,NK_MC_sp,NK_Savescu}.
In this paper we shall
generalize this approach to discretizations of semilinear fractional-parabolic equations.
This, essentially,  will enable us to seamlessly extend the error analysis of the recent paper  \cite{NK_XM}
to the challenging semilinar case and thus obtain
sharp pointwise-in-time error bounds for quasi-graded temporal meshes with
arbitrary degree of grading.
There are a few papers on the numerical analysis of similar nonlinear time-fractional equations
\cite{Du_Y_Zhou,Jin_Li_Zhou,Ji_Liao_L1,Ji_Liao_Al}, but
we are not aware of any such general results in the literature.

The following fractional-in-time semilinear parabolic problem is considered:
\beq\label{problem}
\begin{array}{l}
D_t^{\alpha}u+\LL u+f(x,t,u)=0\quad\mbox{for}\;\;(x,t)\in\Omega\times(0,T],\\[0.2cm]
u(x,t)=0\quad\mbox{for}\;\;(x,t)\in\pt\Omega\times(0,T],\qquad
u(x,0)=u_0(x)\quad\mbox{for}\;\;x\in\Omega.
\end{array}
\eeq
This problem is posed in a bounded Lipschitz domain  $\Omega\subset\R^d$ (where $d\in\{1,2,3\}$).
The operator $D_t^\alpha$, for some $\alpha\in(0,1)$, is
the Caputo fractional derivative in time defined \cite{Diet10} by
\begin{equation}\label{CaputoEquiv}
D_t^{\alpha} u(\cdot,t) :=  \frac1{\Gamma(1-\alpha)} \int_{0}^t(t-s)^{-\alpha}\, \pt_s u(\cdot, s)\, ds
    \qquad\text{for }\ 0<t \le T,
\end{equation}
where $\Gamma(\cdot)$ is the Gamma function, and $\pt_s$ denotes the partial derivative in $s$.
The spatial operator $\LL$ here is a linear second-order elliptic operator:
\beq\label{LL_def}
\LL u := \sum_{k=1}^d \Bigl\{-\pt_{x_k}\!(a_k(x,t)\,\pt_{x_k}\!u) + b_k(x,t)\, \pt_{x_k}\!u \Bigr\}+c(x,t)\, u,
\eeq
with sufficiently smooth coefficients $\{a_k\}$, $\{b_k\}$ and $c$ in $\bar\Omega$, for which we assume that $a_k>0$ in $\bar\Omega$,
and also both $c\ge 0$ and $c-\frac12\sum_{k=1}^d\pt_{x_k}\!b_k\ge 0$.

This problem will be considered under the following assumptions on $f$.

\begin{itemize}
\item[{\bf A1.}]
Let $f$ be continuous in $s$ and  satisfy $f(\cdot,t,s)\in L_\infty(\Omega)$ for all $t>0$ and $s\in\R$,
and the one-sided Lipschitz condition
$$
f(x,t,s_1)-f(x,t,s_2)\ge -\lambda[s_1-s_2]\qquad \forall s_1\ge s_2,\;\; x\in\Omega,\;\;t>0
$$
with some constant $\lambda\ge0$.
%
%
\smallskip

\item[{\bf A2.}]
There exist constants $\sigma_1\le 0\le \sigma_2$  such that $f(\cdot,\cdot,\sigma_1)\le0$ and $f(\cdot,\cdot,\sigma_2)\ge0$,
{\color{red}while $c=0$ in \eqref{LL_def}}.
\end{itemize}
\medskip

Example 1 (Negative {\color{red}reaction} coefficient). The linear  $f=c^*(x,t)\,u+F(x,t)$, with a possibly negative {\color{red}reaction} coefficient $c^*\ge-\lambda$,
clearly satisfies A1.
\smallskip

Example 2 (Allen-Cahn equation). The cubic $f=u^3-u$ satisfies both A1 and A2 with, e.g., $-\sigma_1=\sigma_2=1$.
In particular, the recent papers \cite{Du_Y_Zhou,Ji_Liao_L1,Ji_Liao_Al} are devoted to this equation.
Note that if $|u_0|\le1$, then
$|u|\le1$ $\forall\, t$ \cite[Theorem~2.4]{Du_Y_Zhou}, while our results below imply a similar property for the computed solutions.

\smallskip

Example 3 (Fisher equation). The quadratic $f=u^2-u$  satisfies A2 with, e.g., $\sigma_1=0$ and $\sigma_2=1$, but not A1.
(To be more precise, A2 is satisfied for $s\ge -C$, where $C\ge 0$ is a fixed positive constant.)
Such equations are addressed in \S\ref{ssub_gen_A1}.
\medskip

In this paper, we shall focus
on popular
L1-type schemes for problem \eqref{problem}. Thus,
 consider the discetization of  $D^\alpha_tu$ defined, for $m=1,\ldots,M$, by
\vspace{-0.1cm}
\beq\label{delta_def}
\delta_t^{\alpha} U^m :=  \frac1{\Gamma(1-\alpha)} \sum_{j=1}^m \delta_t U^j\!\int_{t_{j-1}}^{t_j}\!\!(t_m-s)^{-\alpha}\, ds,
\qquad
\delta_t U^j:=\frac{U^j-U^{j-1}}{t_j-t_{j-1}},\vspace{-0.1cm}%
\eeq
when associated with the temporal mesh $0=t_0<t_1<\ldots <t_M=T$ on $[0,T]$.
(Note that, similarly to \cite{NK_L2,NK_XM}, the approach of the present paper may be extended to other discretizations
that are monotone in time.)

 An essential building block in our analysis is
the following  stability result. Given
$\color{red}\lambda\ge 0$ and
$\gamma\in \R$ (where $\gamma\neq 0$ if $\color{red}\lambda>0$), as well as
 a temporal mesh $\{t_j\}_{j=0}^M$ on $[0,T]$ with $\tau:=t_1$,    under certain conditions on the mesh, the following is true for $\{V^j\}_{j=0}^M$:
\beq\label{main_stab}
\left.\begin{array}{c}
|(\delta_t^\alpha-\lambda) V^j|\lesssim (\tau/ t_j)^{\gamma+1}
\\[0.2cm]
\forall\, j\ge1,\;\;\; V^0=0
\end{array}\right\}
\;\;\Rightarrow\;\;
\begin{array}{l}
|V^j|\lesssim {\mathcal V}_\gamma^j:=
\tau t_j^{\alpha-1}\left\{\begin{array}{ll}
1&\mbox{if~}\gamma>0\\
1+\ln(t_j/\tau)&\mbox{if~}\gamma=0\\
(\tau/t_j)^\gamma
&\mbox{if~}\gamma<0
\end{array}\right.\\[-0.3cm]
\;\forall\, j\ge1.\vspace{-0.1cm}
\end{array}
\eeq
The immediate usefulness of this property is due to the fact that truncation errors in time are typically bounded by negative powers of $t_j$.
 Note that
 \eqref{main_stab} is sharp in the sense that it
is consistent with the analogous property for the continuous operator $D_t^\alpha-\lambda$ (similarly to \cite[Remark~1.1]{NK_XM}).
 It is worth mentioning that for $\lambda=0$ it is obtained in \cite{NK_XM} using barrier functions, while here we extend \eqref{main_stab} to $\color{red}\lambda>0$ simply
 as a clever corollary of this property for $\lambda=0$ (by constructing an appropriate upper solution for the operator
 $\delta_t^\alpha-\lambda$).

It should be noted that while the explicit inverse of $D_t^\alpha-\lambda$ is easily available, the proof of \eqref{main_stab} for any discrete operator is quite non-trivial.
As an alternative, discrete Gr\"{o}nwall inequalities were recently employed in the error analysis of L1- and Alikhanov-type schemes \cite{sinum18_liao_et_al,sinum19_liao_et_al,Ji_Liao_L1,Ji_Liao_Al}.
However, the latter approach involves intricate evaluations and, furthermore, yields less sharp error bounds (see Remark~\ref{rem_positive_time} for a more detailed discussion).
Our approach in \cite{NK_XM} and here is entirely different and is substantially more concise as we obtain \eqref{main_stab} essentially using clever barrier functions for $\delta_t^\alpha$,
while the numerical results in \cite{NK_XM} and \S\ref{sec_Num} indicate that our error bounds are sharp in the pointwise-in-time sense.

Similarly to
\cite{Brunner_MC85,Brunner_book,Ji_Liao_L1,Ji_Liao_Al,%
NK_MC_L1,NK_XM,sinum18_liao_et_al,%
mclean_etal_NM07,Mustapha_etal_SINUM14,stynes_etal_sinum17}%
, our main interest will be in graded temporal meshes as they offer an efficient way of computing reliable numerical approximations of solutions singular at $t=0$,
which is typical for \eqref{problem}.
At the same time, as a particular case,
our results immediately apply to uniform temporal meshes.

A number of outstanding theoretical gaps in the error analysis for semilinear fractional-parabolic equations
will be addressed.

\begin{itemize}[leftmargin=0.7cm]

\item
Under very general conditions A1 and A2, whenever the exact solution lies within a certain range (e.g.,
$[\sigma_1,\sigma_2]$, or it is positive),
the  method of discrete upper and lower solutions will easily yield
  a similar property for the computed solutions.
Similar results have been obtained only for the Allen-Cahn equation using the specific form of $f$;
see \cite[Theorem~3.3]{Du_Y_Zhou}
\cite[Theorem~2.2]{Ji_Liao_L1}, \cite[Theorem~3.1]{Ji_Liao_Al}.

\item
Combining the theory of upper and lower solutions with
the subtle and sharp stability property~\eqref{main_stab} will yields sharp pointwise-in-time error bounds for quasi-graded temporal meshes with
arbitrary degree of grading.
We are not aware of any such general results in the literature.

\item
A straightforward particular case of our error bounds is that the (quasi-)uniform temporal mesh yields the first-order convergence in positive time $t\gtrsim 1$ (see Remark~\ref{rem_positive_time}).
This is consistent with the error bounds  in \cite{gracia_etal_cmame,laz_L1,NK_MC_L1,NK_XM} obtained for the linear case,
but appears a new result for the semilinear equations.

\item
Another particular case of our error bounds indicates that
the optimal convergence rates  of order $2-\alpha$ in positive time $t\gtrsim 1$ are attained using much milder
(compared to the optimal) grading
with $r>2-\alpha$
(see Remark~\ref{rem_positive_time}). This is consistent with \cite{NK_XM}, but has not been proved before for the semilinear case.

\item
Note also that
when the optimal grading parameter $r=(2-\alpha)/\alpha$ is used, as particular cases, we recover the optimal global
convergence rate of order $2-\alpha$ (similarly to \cite[Theorem~3.1]{Ji_Liao_L1}),
while in the case of quasi-uniform temporal meshes we recover the global convergence rate of order $\alpha$
(similarly to \cite[Theorem~4.2]{Du_Y_Zhou} and \cite[Theorem~4.4]{Jin_Li_Zhou}); see Remark~
\ref{rem_global_time}.

\end{itemize}

\smallskip

Strictly speaking, Remarks~\ref{rem_positive_time} and~\ref{rem_global_time}, to which we have referred above, apply to the L1 discretizations of the initial-value problem of type~\eqref{problem}. At the same time, the discussion there focuses on the term ${\mathcal E}^m$, which
also appears in the error estimates for
semi-discretizations of the initial-boundary-value problem  \eqref{problem}, 
 and its full discretizations using finite differences 
 and finite elements 
 (see Theorems~\ref{theo_semidiscr}, \ref{theo_FD}, \ref{theo_lumped_FE}, \ref{theo_FE_L2norm} and \ref{theo_FD_Robin}).

 To be more precise with regard to the earlier literature, \cite{Du_Y_Zhou,Ji_Liao_L1,Ji_Liao_Al} are devoted to the Allen-Cahn equation, while \cite{Jin_Li_Zhou} addresses a more general semilinear equation with a Lipschitz-continuous $f=f(u)$ (which is more restrictive compared to A1).
 In \cite{Du_Y_Zhou}, the error is estimated globally in time in the $L_2(\Omega)$ norm for Gr\"{u}nwald-Letnikov-type
 semidiscretizations on uniform temporal meshes.
 In \cite{Jin_Li_Zhou}, similar error bounds are given for  the L1  and the backward Euler convolution quadrature
discretizations in time combined with linear finite elements in space.
In \cite{Ji_Liao_L1,Ji_Liao_Al},  the error is estimated in the $L_\infty(\Omega)$ norm 
for, respectively, the L1 and Alikhanov schemes in time combined with standard finite differences in space
for the case of periodic boundary conditions.




Throughout the paper, it is assumed that 
there exists a unique solution of \eqref{problem} 
such that
{\color{red}at least $u(\cdot,t)\in H^1_0(\Omega)$ $\forall t>0$,}
and
$\|\pt_t^l u(\cdot,t)\|_{L_p(\Omega)}\lesssim 1+t^{\alpha-l}$ for $l\le 2$ with $p\in\{2,\infty\}$.
The latter is a realistic assumption
{\color{red}(e.g., proved in \cite[(3.1)]{Jin_Li_Zhou} for $l=1$, $p=2$)}, 
in contrast to stronger assumptions of type \mbox{$\|{\color{red}\pt^l_t} u(\cdot,t)\|_{L_p(\Omega)}\lesssim 1$} frequently made in the literature.
Indeed,
 \cite[Theorem~2.1]{stynes_too_much_reg} clearly shows
that
this stronger assumption is too restrictive even in the linear case.
{\color{blue}When  full discretizations are considered in \S\S\ref{sec_FD}--\ref{sec_FE}, additional assumptions
are required on  
$\|\pt_{x_k}^l u(\cdot,t)\|_{L_\infty(\Omega)}$ for $l=3,4$ in Theorem~\ref{theo_FD}
and on
$\|\pt_t^l u(\cdot, t)\|_{W^{2}_p(\Omega)}$ for $l=0,1$, $p\in\{2,\infty\}$ in Remarks~\ref{rem_ritz} and~\ref{rem_ritz2}.
For bounds of this type in the linear case, see
\cite[(1.6) and (1.7)]{McL10}, \cite{sakamoto}, 
\cite[\S2]{stynes_etal_sinum17}, \cite[\S6]{NK_MC_L1}.}
%
For some existence, uniqueness and regularity results for the semilinear case, we also refer the reader to
\cite[Theorem~2.3]{Du_Y_Zhou} and \cite[Theorem~3.1]{Jin_Li_Zhou}.
\smallskip

{\it Outline.}
We start by describing
discrete upper and lower solutions and their properties  in \S\ref{sec_upper};
all these results are valid for arbitrarily large $T$ (see Remark~\ref{rem_T_infty}).
Next, \S\ref{sec_L1_stab} is devoted to the proof of the stability result~\eqref{main_stab}
(a version of which for arbitrarily large $T$ is discussed in \S\ref{ssec_long_term}).
This result is then employed  to obtain
pointwise-in-time error bounds for
L1-type discretizations of the initial-value problem of type \eqref{problem} in
 \S\ref{ssec_L1_prdgm},  semi-discretizations of the initial-boundary-value problem  \eqref{problem} in \S\ref{ssec_L1_semi},
 and its full discretizations using finite differences in \S\ref{sec_FD} and finite elements in \S\ref{sec_FE}
 (where the consideration is restricted to $\LL=-\triangle$).
Generalizations of the above results, such as the treatment of other types of boundary conditions, are discussed in \S\ref{sub_gen}.
Finally, our theoretical findings are illustrated by numerical experiments in \S\ref{sec_Num}.
\smallskip

{\it Notation.}
We write
 $a\sim b$ when $a \lesssim b$ and $a \gtrsim b$, and
$a \lesssim b$ when $a \le Cb$ with a generic constant $C$ depending on $\Omega$, $T$, $u_0$,
$f$, and $\alpha$,
but not 
%
 on the total numbers of degrees of freedom in space or time.
  Also, for 
  $1 \le p \le \infty$, and $k \ge 0$,
  we shall use the standard norms
  in the spaces $L_p(\Omega)$ and the related Sobolev spaces $W_p^k(\Omega)$,
  while
  $H^1(\Omega)=W_2^1(\Omega)$ and
  $H^1_0(\Omega)$ is the 
  space of functions in $H^1(\Omega)$ vanishing on $\pt\Omega$.

\section{Discrete upper and lower solutions}\label{sec_upper}
 In this section we shall consider definitions and certain properties of discrete upper and lower solutions
 in the context of the semidiscretization of our original problem.
 Extensions for the operator
 $\delta_t^\alpha-\lambda$ without spatial derivatives
 and certain full discretizations will be given in
 \S\S\ref{ssec2_1}--\ref{ssec2_2}.

Consider the semidiscretization of our problem~\eqref{problem} in time:
\beq\label{semi_semidiscr_method}
\delta_t^\alpha U^m +\LL U^m+ f(\cdot,t_m,U^m)=0\;\;\mbox{in}\;\Omega,\quad U^m=0\;\;\mbox{on}\;\pt\Omega\quad\forall\,m=1,\ldots,M;\quad U^0=u_0.
\eeq

{\sc Definition.} The discrete function $\{\bar U^j\}_{j=0}^M$ is called an {\it upper solution} of problem \eqref{semi_semidiscr_method}
if it satisfies (possibly in a weak sense \cite[\S9.3]{evans}) the following conditions:
\beq\label{upper_semidiscr_method}
\delta_t^\alpha\bar  U^m +\LL\bar  U^m+ f(\cdot,t_m,\bar U^m)\ge0\;\;\mbox{in}\;\Omega,\quad \bar U^m\ge0\;\;\mbox{on}\;\pt\Omega\quad\forall\,m\ge1;\quad \bar U^0\ge u_0.
\eeq
The discrete function $\{\underline{U}^j\}_{j=0}^M$ is called a {\it lower solution} of problem \eqref{semi_semidiscr_method}
if it satisfies the reversed inequalities in \eqref{upper_semidiscr_method}.
\smallskip

\begin{lemma}\label{upper_lemma}
Suppose that $f$ satisfies {\rm A1},
and $\lambda\tau_j^{\alpha}\le \{\Gamma(2-\alpha)\}^{-1}$ $\forall\,j\ge1$.

(i)
If $u_0\in L_\infty(\Omega)$, then problem
\eqref{semi_semidiscr_method} has a unique solution $\{U^j\}_{j=0}^M$, with $U^j\in H^1_0(\Omega)\cap L_\infty(\Omega)$ $\forall j\ge 1$.

(ii) If $\{\bar U^j\}_{j=0}^M$, is an upper solution of problem \eqref{semi_semidiscr_method},
$\bar U^0\in L_\infty(\Omega)$, and $\bar U^j\in H^1(\Omega)\cap L_\infty(\Omega)$ $\forall j\ge 1$,
then
$ U^j\le \bar U^j$ $\forall j\ge 0$.

(iii) If $\{\underline{U}^j\}_{j=0}^M$ is a lower solution of problem \eqref{semi_semidiscr_method},
$\underline{U}^0\in L_\infty(\Omega)$, and $\underline{U}^j\in H^1(\Omega)\cap L_\infty(\Omega)$ $\forall j\ge 1$,
then
$ \underline{U}^j\le  U^j$ $\forall j\ge 0$.
%
%
\end{lemma}%
\smallskip

\begin{proof}
A straightforward calculation shows that
\eqref{delta_def} can be represented as
\beq\label{delta_def_kappa}
\delta^\alpha_tU^m=\kappa_{m,m} U^m-\sum_{j=0}^{m-1}\kappa_{m,j}U^j,
\quad\mbox{where}\quad
\kappa_{m,m}=\frac{\tau_m^{-\alpha}}{\Gamma(2-\alpha)}, \quad \kappa_{m,j}>0\;\;\forall\,m\ge j.
\eeq

(i)
The proof is by induction. Assume that there exist desired $\{U^j\}_{j<m}$.
Combining \eqref{semi_semidiscr_method} with
\eqref{delta_def_kappa},
one concludes that each $U^m$ solves
the semilinear elliptic equation
\beq\label{semi_elliptic}
\N^m U^m:=\LL U^m+\bigl[f(\cdot,t_m, U^m)+\kappa_{m,m} U^m\bigr] = F^m\qquad\mbox{in}\;\Omega,
\eeq
where $F^m:=\sum_{j=0}^{m-1}\kappa_{m,j}U^j$ is a linear combination of  $\{U^j\}_{j<m}$, so
$F^m\in L_\infty(\Omega)$.
As $\lambda\tau_m^{\alpha}\le \{\Gamma(2-\alpha)\}^{-1}$ is equivalent to $\lambda\le \kappa_{m,m}$,
the part $[f(\cdot,\cdot, U^m)+\kappa_{m,m} U^m]$ in \eqref{semi_elliptic} is monotone with respect to $U^m$.
Our assumptions on the elliptic operator $\LL$ imply that it satisfies the maximum principle, so an application of the argument used in the proof of
\cite[Lemma~1]{DK16} with \cite[Lemma 16]{BS73} (where a more general $\LL$ is considered)
to this elliptic equation, subject to $U^m=0$ on $\pt\Omega$,
yields existence of a unique solution $U^m\in H^1_0(\Omega)\cap L_\infty(\Omega)$.
{\color{red}To be more precise, the argument  in \cite[Lemma~1]{DK16}
relies on $\|v\|_{L_\infty(\Omega)}\lesssim \|\LL v\|_{L_\infty(\Omega)}$ $\forall\,v\in H^1_0(\Omega)\cap L_\infty(\Omega)$
 (which follows from the maximum principle \cite[Theorem~3.7 and \S8.1]{GTru}) and essentially
 reduces
 \eqref{semi_elliptic} under assumption A1 (and, hence, with a monotone nonlinearity) to the case addressed in
 \cite[Lemma 16]{BS73}. The latter lemma yields existence of a solution in $H^1_0(\Omega)\cap L_\infty(\Omega)$
for the equation $\LL U^m+g(x,U^m)=0$ with an appropriate $g(x,s)$ uniformly bounded in $\Omega\times\R$, measurable in $x$ and
 continuous in $s$.}
\smallskip

(ii) The proof is again by induction. Assume that we have established  $U^j\le \bar U^j$ for $j<m$.
Then for $\bar U^m$ one gets a version of \eqref{semi_elliptic}:
$\N^m \bar U^m\ge \bar F^m:=\sum_{j=0}^{m-1}\kappa_{m,j}\bar U^j$.
Note that $\bar F^m\ge F^m$ (in view of $\kappa_{m,j}>0$),
so $\N^m \bar U^m\ge \N^m  U^m$.
From this in the domain $\widehat \Omega:=\{U^m>\bar U^m\}$, one gets
$$
\N^m  U^m-\N^m \bar U^m\ge \LL[U^m-\bar U^m]+\underbrace{(\kappa_{m,m}-\lambda)}_{\ge 0}\underbrace{[U^m-\bar U^m]}_{> 0}\ge \LL[U^m-\bar U^m].
$$
Hence
$\LL[U^m-\bar U^m]\le 0$ in $\widehat \Omega$. Finally, an application of the maximum principle for functions in $H^1(\widehat\Omega)$
\cite[\S8.1]{GTru}, one concludes that
$\sup_{\widehat \Omega}(U^m-\bar U^m)\le 0$. The desired bound $U^m\le \bar U^m$ in $\Omega$ follows.
\smallskip

(iii) Imitate the argument of part (ii).
%
%
%
\end{proof}
\smallskip

\begin{corollary}[Bounds for the computed solution]\label{sigma_bounds_U}
Under the conditions of Lemma~\ref{upper_lemma},
suppose that
 $f$ also satisfies {\rm A2}, and
 $\sigma_1\le u_0\le \sigma_2$. Then
 for the unique solution of \eqref{semi_semidiscr_method} one has $\sigma_1\le U^j\le \sigma_2$ $\forall j\ge0$.
\end{corollary}
\smallskip

\begin{proof}
A2 implies that $\sigma_1$ and $\sigma_2$ are, respectively, lower and upper solutions of \eqref{semi_semidiscr_method}. Hence, Lemma~\ref{upper_lemma}(ii),(iii)
yields the desired assertion.
\end{proof}
\smallskip

\subsection{Extension to the operator $\delta_t^\alpha-\lambda$}\label{ssec2_1}~\smallskip

\begin{remark}
\label{rem_upper}
The above definitions of upper and lower solutions, as well as a version of Lemma~\ref{upper_lemma}, but
under a stronger assumption $\lambda\tau_j^{\alpha}<\{\Gamma(2-\alpha)\}^{-1}$ $\forall\,j\ge1$,
clearly apply to
the simpler operator $\delta_t^\alpha-\lambda$ without spatial derivatives.%
\end{remark}%
\smallskip

\begin{corollary}[Comparision principle for $\delta_t^\alpha-\lambda$]\label{cor_comparison}
Let the temporal mesh satisfy $\lambda\tau_j^{\alpha}<\{\Gamma(2-\alpha)\}^{-1}$ $\forall\,j\ge1$.
Then $V^0\le B^0$ and  $(\delta_t^\alpha-\lambda)V^m\le (\delta_t^\alpha-\lambda)B^m$ $\forall\,m\ge1$
imply $V^m\le B^m$ $\forall\,m\ge0$.
\end{corollary}
\smallskip

\begin{proof}
In view of Remark~\ref{rem_upper}, the desired conclusion follows from a version of  Lemma~\ref{upper_lemma}(ii) for the operator $\delta_t^\alpha-\lambda$.
\end{proof}
\smallskip

\subsection{Extension to full discretizations}\label{ssec2_2}
Let
$\bar\Omega_h$ be a finite-dimensional set of points in $\bar\Omega$, comprising the nodes of a certain spatial mesh,
 and $\Omega_h:=\bar\Omega_h\backslash\pt\Omega$ denote the set of interior mesh nodes.
 Consider a fully discrete version of \eqref{semi_semidiscr_method} in the form
\beq\label{FD_problem}
\begin{array}{l}
\delta_t^\alpha U^m(z) +\LL_h U^m(z)+ f(z,t_m,U^m(z))=0\quad\mbox{for}\;\;z\in\Omega_h,\;\;m=1,\ldots,M,\\[0.2cm]
U^m=0\quad\mbox{in}\;\;\bar\Omega_h\cap\pt\Omega,\;\;m=1,\ldots,M,\qquad\;\;
U^0=u_0\quad\mbox{in}\;\;\bar\Omega_h.
\end{array}
\eeq
Generalizing, in an obvious manner, the above definitions  of upper and lower solutions to fully discrete problem~\eqref{FD_problem},
we formulate a version of Lemma~\ref{upper_lemma}.
\smallskip

\begin{lemma}\label{upper_lemma_FD}
Suppose the spatial operator $\LL_h$ in \eqref{FD_problem}
is associated with an M-matrix,
 $f$ satisfies {\rm A1},
and $\lambda\tau_j^{\alpha}\le \{\Gamma(2-\alpha)\}^{-1}$ $\forall\,j\ge1$.

(i)
If $u_0\in L_\infty(\Omega)$, then problem
\eqref{FD_problem} has a unique solution $\{U^j\}_{j=0}^M$.

(ii) If $\{\bar U^j\}_{j=0}^M$, is an upper solution of problem \eqref{FD_problem},
then
$ U^j\le \bar U^j$ $\forall j\ge 0$.

(iii) If $\{\underline{U}^j\}_{j=0}^M$ is a lower solution of problem \eqref{FD_problem},
then
$ \underline{U}^j\le  U^j$ $\forall j\ge 0$.
\end{lemma}
\smallskip

\begin{proof}
For part (i), we imitate the proof of Lemma~\ref{upper_lemma}(i).
For any $m\ge 1$, the solution $U^m$ of \eqref{FD_problem} satisfies
the following version of \eqref{semi_elliptic}:
$\N_h^m U^m:=\LL_h U^m+[f(\cdot,t_m, U^m)+\kappa_{m,m} U^m] = F^m$ in $\Omega_h$.
This is a system of ${\rm dim}(\Omega_h)$ nonlinear equations, and, in view of condition A1 on $f$, the  part $[f(\cdot,\cdot, U^m)+\kappa_{m,m} U^m]$  is monotone in $U^m$.
Consequently, the mapping $\N_h^m$ satisfies the conditions in \cite[\S13.5.6]{orgega_r_book},  which yields existence of a unique solution of this equation $U^m$ in $\Omega_h$.

For parts (ii) and (iii), we start by imitating the proof of Lemma~\ref{upper_lemma}(ii) and, assuming that $U^j\le \bar U^j$ in $\Omega_h$ for $j<m$, conclude that $\N_h^m \bar U^m\ge \N_h^m  U^m$ in $\Omega_h$.
In view of \cite[\S13.5.6]{orgega_r_book}, the mapping $\N_h^m$ is inverse isotone, which immediately  yields
$U^m\le \bar U^m$ in $\Omega_h$.
\end{proof}
\smallskip

\begin{remark}[$T\gg1$]\label{rem_T_infty}
An inspection of the above proofs shows that
all results of \S\ref{sec_upper}, as well as Theorems~\ref{theo_FD}(ii) and~\ref{theo_lumped_FE}(ii), are valid for arbitrarily large $T$.
Additionally, a version of the stability property \eqref{main_stab} for arbitrarily large $T$ will be discussed in \S\ref{ssec_long_term}.%
\end{remark}%

\section{Stability properties of the L1 discrete fractional-derivative operator}\label{sec_L1_stab}

\subsection{Quasi-graded temporal meshes. Main stability result for $\delta^\alpha_t-\lambda$}
Throughout the paper, we shall assume that the temporal mesh is quasi-graded
in the sense that, with some $r\ge 1$,
\beq\label{t_grid_gen}
\tau := t_1\simeq M^{-r},\qquad
\tau_j:=t_j-t_{j-1}
\lesssim\tau^{1/r}t_j^{1-1/r}
\qquad 
\forall\,j=1,\ldots,M.
\eeq
Importantly, the results from \cite{NK_XM}, which we shall employ, apply to this mesh in view of
\cite[Lemma~2.7]{NK_XM}.

For example,  the standard graded temporal mesh
$\{t_j=T(j/M)^r\}_{j=0}^M$ with some $r\ge 1$ (while $r=1$ generates a uniform mesh)
 satisfies \eqref{t_grid_gen},
in view of $\tau_j\simeq M^{-1}\,t_{j-1}^{1-1/r}$ 
and $t_j\le 2^{r}  t_{j-1}$
for $j\ge 2$.

The key  in our error analysis is the following stability property,
which is also the main result of this section.
\smallskip

\begin{theorem}[Stability of $\delta^\alpha_t-\lambda$]\label{theo_new_barrier}
Let  $\lambda\tau_j^{\alpha}<\{\Gamma(2-\alpha)\}^{-1}$ $\forall\,j\ge1$.

(i) Additionally, let the temporal mesh satisfy
 \eqref{t_grid_gen}
with $1\le r\le (2-\alpha)/\alpha$.
Given $\{V^j\}_{j=0}^M$,
the stability property \eqref{main_stab} holds true for any fixed $\lambda\ge 0$ and $\gamma\neq 0$.

(ii)
If $\gamma\le \alpha-1$, then 
one has the above result without assuming \eqref{t_grid_gen}.

(iii)  The above results remain valid if   $|(\delta_t^\alpha-\lambda) V^j|\lesssim (\tau/ t_j)^{\gamma+1}$ in \eqref{main_stab} is replaced by
$(\delta_t^\alpha-\lambda) |V^j|\lesssim (\tau/ t_j)^{\gamma+1}$.
\end{theorem}
\smallskip

Note that the above result is a generalization of the following particular case, addressed in \cite{NK_XM}.
\smallskip
\renewcommand{\thetheorema}{\ref{theo_new_barrier}${}^*$}
\begin{theorema}[
{\cite[Theorem~2.1]{NK_XM}}]\label{theo_barrier}
If $\lambda=0$, then Theorem~\ref{theo_new_barrier} holds true for any fixed $\gamma\in\R$.
\end{theorema}
\smallskip

\subsection{Proof of Theorem~\ref{theo_new_barrier}}
To prove  Theorem~\ref{theo_new_barrier}, we shall employ its particular case, Theorem~\ref{theo_new_barrier}${}^*$ already established in \cite{NK_XM},
and the following lemma.
\smallskip

\begin{lemma}\label{lem_axu_barrier}
For any fixed positive constant $c_0<\frac12\{\lambda \Gamma(2-\alpha)\}^{-1/\alpha}$ such that $\bar\tau:=\max\{\tau_j\}\le\frac12 c_0$,
and any fixed mesh point $t_m\in\{t_j\}_{j=0}^M$, there exists
$\{B^j\}_{j=0}^M$ such that
$$
B^j=0\;\;\forall\,j\le m,\quad
0\le B^j\lesssim 1\quad\mbox{and}\quad
(\delta_t^\alpha- \CCC)\, B^j \gtrsim \left\{\begin{array}{ll}
0&\mbox{for~}t_j< t_m+c_0\\
1&\mbox{for~}t_j\ge t_m+c_0
\end{array}\right.
\;\forall\,j\ge1.
$$
\end{lemma}

Next, we proceed to the proof of Theorem~\ref{theo_new_barrier}, which will be followed by the proof of Lemma~\ref{lem_axu_barrier}.
\smallskip

{\it Proof of Theorem~\ref{theo_new_barrier}.}
(i)
In view of the comparison principle given by Corollary~\ref{cor_comparison}, it suffices to show that
under the conditions of Theorem~\ref{theo_new_barrier}(i)/(ii),
there exists a function $\{\W^j\}$ such that
\beq\label{main_stab_negative_C}
\left.\begin{array}{c}
(\delta_t^\alpha -\CCC){\W}^j\gtrsim (\tau/ t_j)^{1+\gamma}
\\[0.2cm]
\forall j\ge1,\;\;\; {\W}^0=0
\end{array}\right\}
\;\;\mbox{and}\;\;
0\le{\W}^j\lesssim {\mathcal V}_\gamma^j=
\tau t_j^{\alpha-1}
(\tau/t_j)^{\min\{0,\,\gamma\}}
\;\;\forall j\ge 1.
\eeq
Note that here the representation of ${\mathcal V}_\gamma^j$, defined in \eqref{main_stab}, relies on $\gamma\neq 0$.

For any $\gamma\neq 0$, Theorem~\ref{theo_barrier}(i), the conditions of which are also satisfied, yields
\beq\label{main_stab_B}
\left.\begin{array}{c}
\delta_t^\alpha {\B}_\gamma^j= (\tau/ t_j)^{1+\gamma}
\\[0.2cm]
\forall j\ge1,\;\;\; {\B}_\gamma^0=0
\end{array}\right\}
\;\;\Rightarrow\;\;
0\le{\B}_\gamma^j\lesssim {\mathcal V}_\gamma^j=
\tau^\alpha
(\tau/t_j)^{1+\min\{0,\,\gamma\}-\alpha}
\;\;\;\forall j\ge 1.
\eeq
Here the representation of ${\mathcal V}_\gamma^j$ is different from (but equivalent to) the one in \eqref{main_stab_negative_C},
and will be more convenient in what follows.

Set
$$
\gamma^*:=\min\{0,\,\gamma\}-\alpha<0.
$$
Now, \eqref{main_stab_B} implies that,
for a sufficiently large constant $C$,
\begin{align*}
0\le \B_{\gamma}^j&\lesssim \tau^\alpha(\tau/ t_j)^{1+\gamma^*}&\Rightarrow&&
(\delta_t^\alpha -\CCC)\,{\B}_{\gamma}^j
&\gtrsim(\tau/ t_j)^{1+\gamma}-C\tau^\alpha(\tau/ t_j)^{1+\gamma^*},
\\[0.2cm]
0\le \B_{\gamma^*}^j&\lesssim \tau^\alpha(\tau/ t_j)^{1+\gamma^*-\alpha}&\Rightarrow&&
(\delta_t^\alpha -\CCC)\,{\B}_{\gamma^*}^j
&\gtrsim(\tau/ t_j)^{1+\gamma^*}-C\tau^\alpha(\tau/ t_j)^{1+\gamma^*-\alpha}
\\
&&&&&\gtrsim(\tau/ t_j)^{1+\gamma^*}\bigl[1-C t_j^{\alpha}\bigr],
\end{align*}
where we also used
$\gamma^{**}:=\min\{0,\,\gamma^*\}-\alpha=\gamma^*-\alpha$.
Consequently, for a sufficiently large constant $\bar c$ and a sufficiently small constant $c_1$, one obtains
$$
(\delta_t^\alpha -\CCC)\,\bigl[\B_{\gamma}^j+\bar c\,\tau^\alpha\,{\B}_{\gamma^*}^j\bigr]
\gtrsim(\tau/ t_j)^{1+\gamma}-
C\left\{\begin{array}{cl}
0&\mbox{for~}t_j< c_1,\\[0.2cm]
\tau^{1+\min\{0,\,\gamma\}}&\mbox{for~}t_j\ge c_1.
\end{array}\right.
$$
Here, for the case $t_j\ge c_1$, we also employed $\tau^\alpha\,\tau^{1+\gamma^*}=\tau^{1+\min\{0,\,\gamma\}}$.
Note also that
\beq\label{upper_B_gamma}
\B_{\gamma}^j+\bar c\,\tau^\alpha\,{\B}_{\gamma^*}^j\lesssim \tau^\alpha(\tau/ t_j)^{1+\gamma^*}={\mathcal V}_\gamma^j\,.
\eeq

Finally, let
$$
\W^j:=\B_{\gamma}^j+\bar c\,\tau^\alpha\,{\B}_{\gamma^*}^j+\bar c^2\,\tau^{1+\min\{0,\,\gamma\}}\,B^j,
$$
where $\{B^j\}$ is from Lemma~\ref{lem_axu_barrier} with
$c_0:=\frac12 c_1$ and
$t_m\simeq 1$
such that $t_m+c_0\le c_1$.
Then $(\delta_t^\alpha -\CCC)\W^j\gtrsim(\tau/ t_j)^{1+\gamma}$, in agreement with \eqref{main_stab_negative_C},
while  the required upper bound $\W^j\lesssim {\mathcal V}_\gamma^j$
follows from \eqref{upper_B_gamma} combined with
$B^j=0$ for $j\le m$ and $\tau^{1+\min\{0,\,\gamma\}} B_j\lesssim \tau t_j^{\alpha-1}
(\tau/t_j)^{\min\{0,\,\gamma\}}= {\mathcal V}_\gamma^j$ for $t_j\ge t_m\gtrsim1$.
Thus, \eqref{upper_B_gamma} is established.
\smallskip

(ii)
Note that
as now the conditions of Theorem~\ref{theo_barrier}(ii) are satisfied, one gets~\eqref{main_stab_B} only for any $\gamma\le\alpha-1$.
Importantly,
if $\gamma$ satisfies the latter restriction, so does~$\gamma^*$.
Hence, the proof of part (i) applies to this case.
\smallskip

(iii)
Let $W^0=0$ and $(\delta_t^\alpha-\lambda) W^j
\simeq (\tau/ t_j)^{\gamma+1}
\ge
(\delta_t^\alpha-\lambda) |V^j|$  $\forall\,j\ge1$.
Then $0\le |V^j|\le W^j$ $\forall\,j\ge1$
(in view of Corollary~\ref{cor_comparison}), while  the results of parts (i) and (ii) apply to $\{W^j\}$.
\endproof
\smallskip

It remains to prove the auxiliary Lemma~\ref{lem_axu_barrier} (which we used in the above proof).%
\smallskip

{\it Proof of Lemma~\ref{lem_axu_barrier}.}
 First, consider the case $t_m=0$. Let
 $$
B(t):=\sum_{k=0}^K {\bar c}^k\, B_k(t),\quad
B_k(t):=\max\{0,\, t-q_k\},\quad q_0:=0,\quad q_k\in[c_0 k-\bar\tau,c_0 k].
$$
Here  $0\le K\lesssim 1$ is chosen so that $T\in(q_{K+1},q_{K+2}]$, i.e.
$K+2 = \bigl\lceil T/c_0 \bigr\rceil\lesssim 1$
(unless $c_0\ge T$, in which case $K:=0$).

Applying the continuous operator $D_t^\alpha- \CCC$
to
 $B_0=t$, one easily gets
\begin{align*}
(D_t^\alpha- \CCC)\, B_0(t)=t^{1-\alpha}\bigl(\{ \Gamma(2-\alpha)\}^{-1}-\CCC t^\alpha\bigr)
&\gtrsim \left\{\begin{array}{rl}
0&\mbox{for~}t\in(0,c_0)\supset(0,q_1),\\
1&\mbox{for~}t\in [c_0-\bar\tau,2c_0]\supset[q_1,q_2],\\
\!\!\!-1&\mbox{for~}t\in(q_2,T].
\end{array}\right.
\intertext{In a similar manner, $\forall\,m$ one gets}
(D_t^\alpha- \CCC)\, B_m(t)
&\gtrsim \left\{\begin{array}{rl}
0&\mbox{for~}t\in(0,q_{m+1}),\\
1&\mbox{for~}t\in [q_{m+1},q_{m+2}],\\
\!\!\!-1&\mbox{for~}t\in(q_{m+2},T].
\end{array}\right.
\intertext{Note that $T\in(q_{K+1},q_{K+2}]$ implies $(q_{m+2},T]=\emptyset$ for $m=K$.
As
$K\lesssim 1$, choosing $\bar c\lesssim 1$ sufficiently large in the definition of $B(t)$, one can obtain}
(D_t^\alpha- \CCC)\,B(t)
&\gtrsim\left\{\begin{array}{rl}
0&\mbox{for~}t\in(0,q_1),\\
1&\mbox{for~}t\in [q_1,T].
\end{array}\right.
\end{align*}
It remains to choose
$\{q_m\}_{m=0}^K\subset\{t_j\}_{j=0}^M$, e.g., by
letting each $q_m$ be the maximal mesh point subject to $q_m\in[c_0 m-\bar\tau,c_0 m]$.
Then $(\delta_t^\alpha- \CCC)B(t_j)=(D_t^\alpha- \CCC)B(t_j)$, and the desired result follows
for the discrete function $B^j:=B(t_j)$.

 Finally, consider the case $t_m>0$.
On the sub-mesh $\{t_j\}_{j=m}^M$, construct a discrete function as above. Augmenting this function by zeros on the remaining sub-mesh $\{t_j\}_{j=0}^{m-1}$,
one gets the desired $\{B^j\}_{j=0}^M$.
\endproof

{\color{blue}
\subsection{Long-time stability of $\delta^\alpha_t-\lambda$ for $\gamma+1=\alpha$}\label{ssec_long_term}

The stability property \eqref{main_stab} is established in Theorem~\ref{theo_new_barrier} under the assumption that $T\lesssim 1$, i.e.
the terminal time is bounded. At the same time, long-time solutions are frequently of interest for semilinear problems.
However,
the analysis of long-time stability and convergence is considerably more challenging even for classical parabolic equations.
Here we shall extend \eqref{main_stab} to
the case of arbitrarily large $T$ and arbitrary temporal meshes for one particular case
of $\gamma+1=\alpha$ (which corresponds to the optimal grading parameter $r=(2-\alpha)/\alpha$ in the convergence analysis of Section~\ref{ssec_L1_prdgm}).

In the remainder of this section, all constants, including those used in the definition of the notation of type $\lesssim$, will be understood as independent of $T$.
\smallskip
\begin{lemma}
Suppose that $\gamma+1=\alpha$, and $T$ is arbitrarily large.\\
(i) If $\lambda=0$, then \eqref{main_stab} remains valid on an arbitrary mesh independently of $T$.\\
(ii)
For any $\lambda'>\lambda>0$, there exist $c_0=c_0(\lambda')$ and $\bar C=\bar C(\lambda')$, independent of $T$ and $m$, such that
if $\tau_j\le c_0$ $\forall j\le m$, then
\beq\label{main_stab_infty}
V^0=0,\;\; |(\delta_t^\alpha-\lambda) V^j|\le (\tau/t_j)^\alpha\;\;\forall\, j\ge1
\;\;\Rightarrow\;\;
|V^j|\le
\bar C\tau^\alpha\, E_\alpha(\lambda' t^\alpha_j)
\;\;\forall\, j\le m,
\eeq
where $E_\alpha(s):=\sum_{k=0}^\infty\frac{s^k}{\Gamma(k\alpha+1)}$ is
the Mittag-Leffler function of order $\alpha$.
\end{lemma}
\smallskip

\begin{proof}
(i)
For $\lambda=0$,   \eqref{main_stab} involves ${\mathcal V}_{\alpha-1}^j=\tau^\alpha$ and remains valid on an arbitrary mesh
for an arbitrarily large $T$, as can be shown by an inspection of the proof of
\cite[Lemma~2.1]{NK_MC_L1}. The latter proof may be interpreted as using the barrier  $\{\mathbbm{1}^j\}_{j\ge0}$ defined by  $\mathbbm{1}^0=0$ and $\mathbbm{1}^j=1$ for $j\ge1$, which, combined with \eqref{delta_def},
yields
$\delta_t^\alpha \mathbbm{1}^j\gtrsim t_j^{-\alpha}$ independently of $T$. The desired bound  follows.

(ii)
Let $\lambda'>\lambda>0$ and, to simplify the presentation, first consider
a uniform mesh $\{t_j=j\tau\}_{j=0}^\infty$.
Clearly, $(\delta_t^\alpha-\lambda) \mathbbm{1}^j\ge C t_j^{-\alpha}-\lambda$.
 Set
 $\mathbbm{E}(t):=E_\alpha(\lambda' t^\alpha)-1$,
for which \cite[Theorem~4.3]{Diet10} yields $(D_t^\alpha -\lambda' )\mathbbm{E}(t)=\lambda'$.
Next, splitting $\sum_{k=1}^{\infty}$ in the definition of $E_\alpha$ leads to
$\mathbbm{E}=\invbreve{\mathbbm{E}}+\breve{\mathbbm{E}}=:\sum_{1\le k\le 1/\alpha} +\sum_{k>1/\alpha}$,
where $\invbreve{\mathbbm{E}}$ is concave and $\breve{\mathbbm{E}}$ is convex.
For a generic function $v$, let
$v^I$ be its piecewise-linear interpolant, and
 $v^*(t):=v(t-\tau)$ for $t\ge \tau$ with  $v^*(t):=0$ otherwise.

We shall show below that if $\tau$ is sufficiently small, then
 \begin{subequations}\label{star_long}
\begin{align}
\label{star_long_a}
&D_t^\alpha(\invbreve{\mathbbm{E}}^*+\breve{\mathbbm{E}})^I
\ge
D_t^\alpha(\invbreve{\mathbbm{E}}^*+\breve{\mathbbm{E}}^*) =D_t^\alpha{\mathbbm{E}}^*
=\lambda' ({\mathbbm{E}}^*+1)
\qquad\forall\,t_j,\;j\ge2,
\\
\label{star_long_b}
&
\lambda' ({\mathbbm{E}}^*+1)
\ge
\lambda ({\mathbbm{E}}+{\textstyle\frac12})
\ge \lambda (\invbreve{\mathbbm{E}}^*+\breve{\mathbbm{E}})
+{\textstyle\frac12}\lambda
\qquad\qquad\qquad\;\;\;\forall\,t>0.
\end{align}
\end{subequations}
Now, for the barrier function $B^j:= \mathbbm{1}^j+2(\invbreve{\mathbbm{E}}^*+\breve{\mathbbm{E}})(t_j)$, the bounds \eqref{star_long} imply that
$(\delta_t^\alpha-\lambda)B^j\ge\delta_t^\alpha\mathbbm{1}^j\gtrsim t_j^{-\alpha}$ $\forall\,j\ge 2$, while
 $B^j\le 1+2\mathbbm{E}(t_j)\le 2E_\alpha(\lambda' t^\alpha_j)$.
Note also that $(\delta_t^\alpha-\lambda)B^1\gtrsim(\tau^{-\alpha}-\lambda)B^1\gtrsim t_1^{-\alpha}$.
The desired bound \eqref{main_stab_infty} follows.

It remains to establish \eqref{star_long}. 
For \eqref{star_long_a}, note that $D_t^\alpha((\invbreve{\mathbbm{E}}^*)^I-\invbreve{\mathbbm{E}}^*)(t_j)\ge 0$ follows from
$\invbreve{\mathbbm{E}}^*-(\invbreve{\mathbbm{E}}^*)^I\ge 0$ (the latter in view of $\invbreve{\mathbbm{E}}$ being concave),
and can be shown by recalling  \eqref{CaputoEquiv} and then applying integration by parts
(similarly to the truncation error representations used in the proofs of \cite[Lemma~2.3]{NK_MC_L1} and \cite[Lemma~3.4]{NK_XM}).
We also use  \mbox{$D_t^\alpha(\breve{\mathbbm{E}}^I-\breve{\mathbbm{E}}^*)(t_j)\ge 0$,} which follows from \eqref{CaputoEquiv} combined with
the observation that
on each $(t_{j-1},t_j)$ one has $\frac{d}{dt}\breve{\mathbbm{E}}^I(t)\ge \frac{d}{dt}\breve{\mathbbm{E}}(t_{j-1})=
\frac{d}{dt}\breve{\mathbbm{E}}^*(t_j)\ge \frac{d}{dt}\breve{\mathbbm{E}}^*(t)$ (as $\breve{\mathbbm{E}}$ is convex).
The final inequality in \eqref{star_long_b} follows from ${\mathbbm{E}}=\invbreve{\mathbbm{E}}+\breve{\mathbbm{E}}$ combined with
$\invbreve{\mathbbm{E}}\ge \invbreve{\mathbbm{E}}^*$.
Finally, for the first bound in \eqref{star_long_b},
note that
the asymptotic representation \cite[(1.8.27)]{MittagL_book} of $E_\alpha(s)$ implies that, with a sufficiently large constant $c_2$,
for $t-\tau\ge c_2$ one has
$|{\mathbbm{E}}(t)+1-\alpha^{-1}\exp(\lambda'^{1/\alpha}t)| \le\frac14\lambda/\lambda'$.
Using the latter to estimate $\lambda' ({\mathbbm{E}}^*+1)-\lambda ({\mathbbm{E}}+1)+\frac12\lambda$,
 yields \eqref{star_long_b} for $t\ge c_2$ provided that $\tau\le \lambda'^{-1/\alpha}\ln (\lambda'/\lambda)$.
Additionally,  ${\mathbbm{E}}(t)+\frac12\le 1\le {\mathbbm{E}}^*(t)+1$ for $t\le c_1$,
with a sufficiently small $c_1$,
again yields \eqref{star_long_b} for $t\le c_1$. For the remaining  $t\in(c_1,c_2)$, using $|{\mathbbm{E}}'|\lesssim 1$
one gets ${\mathbbm{E}}^*-{\mathbbm{E}}\gtrsim-\tau\ge -\frac12 $, so \eqref{star_long_b} is proved $\forall\,t>0$.

If the mesh is nonuniform, a version of the above argument applies with
$v^*(t):=v(t-t_l^*)$, where $t_l^*$ is the minimal node in $\{t_j\}$ such that $t_l^*\ge\bar\tau_m:=\max_{j\le m}\tau_j$ (hence, $t_l^*\le
t_{l-1}+\bar\tau_m\le 2\bar\tau_m$ is sufficiently small).
\end{proof}
\smallskip

\begin{remark}\label{rem_stab_sharp}
If  $\lambda'>\lambda$ in the stability result \eqref{main_stab_infty} is replaced by $\lambda$, then
it
becomes consistent with the analogous property for the continuous Caputo operator $D_t^\alpha$.
Indeed,
for $\mathbbm{1}(t):=\min\{t/\tau,\,1\}$
and $\bar{\mathbbm{E}}(t):=E_\alpha(\lambda t^\alpha)-1$, 
a calculation shows that if $t\ge \tau$, then
$(D_t^\alpha-\lambda)[\mathbbm{1}+\bar{\mathbbm{E}}](t)=D_t^\alpha\mathbbm{1}(t)
\simeq\tau^{-1}[t^{1-\alpha}-(t-\tau)^{1-\alpha}]\simeq t^{-\alpha}$
and $[\mathbbm{1}+\bar{\mathbbm{E}}](t)=E_\alpha(\lambda t^\alpha)$.
\end{remark}

}


\section{Error estimation for a simplest example\! (without spatial derivatives)}\label{ssec_L1_prdgm}
It is convenient to illustrate our approach to the estimation of the temporal-discretization error using a very simple example.
Consider a semilinear fractional-derivative problem without spatial derivatives together with its discretization:
\begin{subequations}\label{simplest}
\begin{align}
\label{simplest_a}
D_t^\alpha u(t)\;+\;\;f(t,u)\;\;&=0&&\hspace{-1.6cm}\mbox{for}\;\;t\in(0,T],&&\hspace{-0.8cm} u(0)=u_0,
\\ \delta_t^\alpha U^m+f(t_m,U^m)&=0&&\hspace{-1.6cm}\mbox{for}\;\;m=1,\ldots,M,&&\hspace{-0.8cm} U^0=u_0.
\label{simplest_b}
\end{align}
\end{subequations}
Throughout this section, with slight abuse of notation, $\pt_t $ will be used for $\frac{d}{dt}$.

The main result here is the following error estimate. 
\smallskip

\begin{theorem}\label{lem_simplest_star}
(i)
Let the temporal mesh satisfy
 \eqref{t_grid_gen}
with $r\ge 1$, and let $\lambda\tau_j^{\alpha}<\{\Gamma(2-\alpha)\}^{-1}$ $\forall\,j\ge1$.
Suppose that $u$ 
is a unique solution of  \eqref{simplest_a},
in which $f$ satisfies a version of {\rm A1},
and
$|\partial_t^l u |\lesssim 1+t^{\alpha-l}$ for $l = 1,2$ and $t\in(0,T]$.
Then
there exists a unique solution $\{U^m\}$ of \eqref{simplest_b}, and
$\forall\,m\ge1$
\beq\label{E_cal_m}
|u(t_m)-U^m|\lesssim {\mathcal E}^m:=
\left\{\begin{array}{ll}
M^{-r}\,t_m^{\alpha-1}&\mbox{if~}1\le r<2-\alpha,\\[0.3cm]
M^{-r(1-\epsilon)}\,t_m^{\alpha-(1-\epsilon)}
&\mbox{if~}r=2-\alpha,\\[0.3cm]
M^{\alpha-2}\,t_m^{\alpha-(2-\alpha)/r}
&\mbox{if~}r>2-\alpha,
\end{array}\right.
\eeq
where $\epsilon$ is an arbitrarily small positive constant.

(ii)
If, additionally, $f$ satisfies a version of {\rm A2},
and $\sigma_1\le u_0\le\sigma_2$,
then $\sigma_1\le U^m\le\sigma_2$ $\forall\, m\ge 0$.
\end{theorem}
\smallskip

\begin{remark}[Case $r=2-\alpha$]
Note that for the case $\lambda =0$ in A1, one can easily get
a slightly sharper version of \eqref{E_cal_m} with
$$
{\mathcal E}^m:=M^{\alpha-2}\,t_m^{\alpha-1}[1+\ln(t_m/t_1)]\qquad \mbox{if~}r=2-\alpha\mbox{~~and~~}\lambda=0
$$
(similarly to the results for the linear case in \cite{NK_XM}).
In comparison, \eqref{E_cal_m} gives a slightly less optimal bound because we have established \eqref{main_stab} for $\gamma=0$ only when $\lambda=0$
(see Theorem~\ref{theo_barrier}).
\end{remark}
\smallskip

\begin{remark}[Convergence in positive time]\label{rem_positive_time}
Consider $t_m\gtrsim 1$. Then
${\mathcal E}^m\simeq M^{-r}$ for $r<2-\alpha$ and ${\mathcal E}^{m}\simeq M^{\alpha-2}$ for $r>2-\alpha$,
i.e. in the latter case the optimal convergence rate is attained.
For $r=2-\alpha$ one gets an almost optimal convergence rate as now ${\mathcal E}^{m}\simeq M^{(\alpha-2)(1-\epsilon)}$ with an arbitrarily small $\epsilon>0$.

Note also that for $r=1$ (i.e. for the quasi-uniform temporal mesh), we have ${\mathcal E}^m\simeq M^{-1}$.
This is consistent with the error bounds  in \cite{gracia_etal_cmame,laz_L1,NK_MC_L1,NK_XM} obtained for the linear case,
but appears a new result for the semilinear equations.

By contrast, \cite[Theorem~3.1]{Ji_Liao_L1}
(obtained by means of a discrete Gr\"{o}nwall inequality for the time-fractional Allen-Cahn equation)
 gives a somewhat similar, but considerably less sharp error bound for graded meshes, as (in our notation)
it involves the term $O(\tau^\alpha)=O(M^{-\alpha r})$, so it requires (in our notation) $r=(2-\alpha)/\alpha$ to attain the optimal convergence rate in positive time.
In fact, for any $ r<(2-\alpha)/\alpha$, our error bound is sharper than the pointwise-in-time bound from \cite[Theorem~3.1]{Ji_Liao_L1}.
(Note also that a similar term $O(\tau^\alpha)=O(M^{-\alpha r})$ appears in the
error bound of
 \cite[ Theorem~4.1]{Ji_Liao_Al} for the higher-order Alikhanov scheme.)
\end{remark}
\smallskip

\begin{remark}[Global convergence]\label{rem_global_time}
Note that $
\max_{m\ge1}{\mathcal E}^m\simeq {\mathcal E}^1\simeq \tau_1^\alpha\simeq M^{-\alpha r}$ for $\alpha\le (2-\alpha)/r$,
while
$\max_{m\ge1}{\mathcal E}^m\simeq {\mathcal E}^M\simeq M^{\alpha-2}$ otherwise.
Consequently, Theorem~\ref{lem_simplest_star} yields the global error bound
$|u(t_m)-U^m|\lesssim M^{-\min \{\alpha r,2-\alpha\}} $.

This immediately implies that
the optimal grading parameter for global accuracy is $r=(2-\alpha)/\alpha$.
Note that similar global error bounds were obtained in \cite{sinum18_liao_et_al,NK_MC_L1,stynes_etal_sinum17} for the linear case,
and in \cite[Theorem~3.1]{Ji_Liao_L1} for the Allen-Cahn equation.

For $r=1$, our global error bound becomes
$|u(t_m)-U^m|\lesssim M^{-\alpha } $, which is consistent with the bounds of
\cite[Theorem~4.2]{Du_Y_Zhou} and
\cite[Theorem~4.4]{Jin_Li_Zhou}, respectively
obtained for Gr\"{u}nwald-Letnikov-type schemes
and
for L1-type schemes, as well as for the backward Euler convolution quadrature.
\end{remark}
\smallskip

{\it Proof of Theorem~\ref{lem_simplest_star}.}
In view of Remark~\ref{rem_upper}, the existence of a unique solution $\{U^m\}$ follows from Lemma~\ref{upper_lemma}(i), while part (ii)
follows from Corollary~\ref{sigma_bounds_U}.

It remains to establish \eqref{E_cal_m}.
Consider the error $e^m:=u(t_m)-U^m$
and the truncation error $r^m:=\delta_{t}^{\alpha} u(t_m)-D_t^\alpha u(t_m)$ $\forall\,m\ge 1$.
A standard calculation using \eqref{simplest} yields
$e^0=0$ and
\beq\label{error_eq_simple}
\delta_t^\alpha e^m+[f(t_m,U^m+e^m)-f(t_m,U^m)]=r^m\qquad \forall\,m\ge 1.
\eeq
Multiply this equation by $\varsigma^m:={\rm sign}(e^m)$
and note that $\varsigma^me^m=|e^m|$ so
\begin{align*}
\varsigma^m(\delta_t^\alpha e^m)&\ge\kappa_{m,m}|e^m|-\sum_{j=0}^{m-1}\underbrace{\kappa_{m,j}}_{>0} |e^{j}|= \delta_t^\alpha |e^m|,
\\
\varsigma^m[f(t_m,U^m+e^m)-f(t_m,U^m)]&\ge -\lambda |e^m|,
\end{align*}
where we used \eqref{delta_def_kappa} and condition A1 on $f$.
Hence,
we arrive at
\beq\label{error_ineq}
(\delta_t^\alpha -\lambda)\, |e^m|\le |r^m|\qquad \forall\,m\ge 1.
\eeq

For the truncation error, recall from \cite[Lemma~3.4 and proof of Theorem~3.1]{NK_XM} that
\beq\label{r_bound_simple}
|r^m|\lesssim (\tau/t_m)^{\gamma+1}\quad\forall\,m\ge 1,
\qquad\mbox{where}\;\;
\gamma+1:=\min\{\alpha+1,(2-\alpha)/r\}.
\eeq
Hence, we can apply Theorem~\ref{theo_new_barrier} to \eqref{error_ineq} (in particular, note part (iii) of this theorem).

Consider three cases.
\smallskip

Case $1\le r<2-\alpha$. Then both $(2-\alpha)/r> 1$ and $\alpha+1>1$,
so $\gamma>0$.
An application of Theorem~\ref{theo_new_barrier}(i) for this case yields
$|e^m|\lesssim \tau\, t_m^{\alpha-1}$, where $\tau\sim M^{-r}$.
\smallskip

Case $r=2-\alpha$. Then  $(2-\alpha)/r= 1$, while $\alpha+1>1$,
so $\gamma=0$. As our stability result does not apply to this case, we note that now
$|r^m|\lesssim \tau/t_m\lesssim (\tau/t_m)^{1-\epsilon}$ for an arbitrarily small $\epsilon>0$.
An application of Theorem~\ref{theo_new_barrier}(i) yields
$|e^m|\lesssim \tau\, t_m^{\alpha-1}(\tau/t_m)^{-\epsilon}\sim
\tau^{1-\epsilon}\,t_m^{\alpha-(1-\epsilon)}
$,
where $\tau\sim M^{-r}$, so $\tau^{1-\epsilon}\sim M^{-r(1-\epsilon)}$

\smallskip

Case $r>2-\alpha$.
Then
$(2-\alpha)/r< 1$, while $\alpha+1>1$,
so $\gamma+1=(2-\alpha)/r<1$.
An application of Theorem~\ref{theo_new_barrier}(where part~(i) of this theorem is used if $r\le (2-\alpha)/\alpha$
and part~(ii) is used otherwise)
 yields
$|e^m|\lesssim \tau\, t_m^{\alpha-1}(\tau/t_m)^{(2-\alpha)/r-1}\sim
\tau^{(2-\alpha)/r}t_m^{\alpha-(2-\alpha)/r}$,
where
$\tau^{(2-\alpha)/r}\sim M^{\alpha-2}$.
\endproof

\section{Error analysis for the L1 semidiscretization in time}\label{ssec_L1_semi}

Recall the semi\-discretization of our problem~\eqref{problem} in time, given by \eqref{semi_semidiscr_method}.
\smallskip

\begin{theorem}\label{theo_semidiscr}
(i)
Let the temporal mesh satisfy
 \eqref{t_grid_gen}
with $r\ge 1$, and let $\lambda\tau_j^{\alpha}<\{\Gamma(2-\alpha)\}^{-1}$ $\forall\,j\ge1$.
Suppose that $u$ is a unique solution of \eqref{problem},\eqref{LL_def} with the initial condition $u_0\in L_\infty(\Omega)$ and
 under assumption {\rm A1} on $f$.
Also,  given $p\in\{2,\infty\}$,
suppose that
$\color{red}u(\cdot,t)\in H^1_0(\Omega)$ for $t\in(0,T]$
and
$\|\partial_t^l u (\cdot, t)\|_{L_p(\Omega)}\lesssim 1+t^{\alpha-l}$ for $\color{red}l = 0,1,2$ and $t\in(0,T]$.
Then there exists a unique solution  $\{U^m\}$ of \eqref{semi_semidiscr_method}, and
\beq\label{L1_semi_error}
\|u(\cdot,t_m)-U^m\|_{L_p(\Omega)}\lesssim {\mathcal E}^m
\qquad\forall\,m=1,\ldots,M, 
\eeq
where ${\mathcal E}^m$ is from \eqref{E_cal_m}.

(ii)
If, additionally, $f$ satisfies {\rm A2},
and $\sigma_1\le u_0\le\sigma_2$,
then $\sigma_1\le U^m\le\sigma_2$ $\forall\, m\ge 0$.
\end{theorem}
\smallskip

\begin{proof}
We imitate the proof of Theorem~\ref{lem_simplest_star}.

The existence of a unique solution $U^m\in H^1_0(\Omega)\cap L_\infty(\Omega)$ for $m\ge 1$ follows from Lemma~\ref{upper_lemma}(i), while part (ii)
follows from Corollary~\ref{sigma_bounds_U}.

It remains to establish~\eqref{L1_semi_error}.
For the error $e^m:= u(\cdot,t_m)-U^m\in \color{red}H^1_0(\Omega)\cap L_p(\Omega)$, using \eqref{problem} and \eqref{semi_semidiscr_method},
one gets
$e^0=0$ and
\beq\label{semidsicr_em_eq}
\delta_t^\alpha e^m +\LL e^m+[f(\cdot, t_m,U^m+e^m)-f(\cdot, t_m,U^m)]=r^m\qquad \forall\,m\ge 1
\eeq
(which is a version of \eqref{error_eq_simple}).
Here 
$r^m:=\delta_{t}^{\alpha} u(\cdot,t_m)-D_t^\alpha u(\cdot, t_m)$, and, similarly to \eqref{r_bound_simple}, it satisfies
\beq\label{trunc_er_semidis}
\|r^m\|_{L_p(\Omega)}\lesssim (\tau/t_m)^{\gamma+1}\quad\forall\,m\ge 1,
\quad\mbox{where}\;\;
\gamma+1:=\min\{\alpha+1,(2-\alpha)/r\}.
\eeq
Hence, to get the desired bound~\eqref{L1_semi_error}
it suffices to prove
\beq\label{semi_error}
(\delta_t^\alpha-\lambda) \|e^m\|_{L_p(\Omega)}\le \|r^m\|_{L_p(\Omega)}
\qquad\forall\,m\ge 1,
\eeq
which is a version of \eqref{error_ineq},
so one then proceeds as in the proof of the error bound \eqref{E_cal_m} in Theorem~\ref{lem_simplest_star}.
 The cases $p=2$ and $p=\infty$ of \eqref{semi_error} will be addressed separately.

For $p=2$, consider the $L_2(\Omega)$ inner product (denoted $\langle\cdot,\cdot\rangle$)
of \eqref{semidsicr_em_eq} with $e^m$.
Clearly,
$c-\frac12\sum_{k=1}^d\pt_{x_k}\!b_k\ge 0$ implies
$\langle\LL e^m,e^m\rangle\ge0$, and we also have $\langle r^m,\,e^m\rangle\le \|r^m\|_{L_2(\Omega)}\|e^m\|_{L_2(\Omega)}$.
Furthermore, recalling \eqref{delta_def_kappa} and condition A1 on $f$, one concludes that
\begin{subequations}\label{semi_aux_p2}\begin{align}
\bigl\langle\delta_t^\alpha e^m,\,e^m\bigr\rangle =\kappa_{m,m}\|e^m\|^2_{L_2(\Omega)}-\sum_{j=0}^{m-1}\underbrace{\kappa_{m,j}}_{>0} \bigl\langle e^{j},\,e^m\bigr\rangle
&\ge  \bigl(\delta_t^\alpha \|e^m\|_{L_2(\Omega)}\bigr)\|e^m\|_{L_2(\Omega)}\,.
\\[0.2cm]
\bigl\langle f(\cdot, t_m,U^m+e^m)-f(\cdot, t_m,U^m),\,e^m\bigr\rangle&\ge-\lambda \|e^m\|^2_{L_2(\Omega)}\,.
\end{align}\end{subequations}
Combining these findings, one gets \eqref{semi_error} for $p=2$.

For $p=\infty$, in view of our assumptions,  \eqref{semidsicr_em_eq} implies that $\LL e^m\in L_{\infty}(\Omega)$
{\color{red}(where the bounds of type $ f(\cdot, t_m, U^m)\le f(\cdot,t_m,\|U^m\|_{L_\infty(\Omega)})+2\lambda\|U^m\|_{L_\infty(\Omega)}$ were used).}
So $e^m\in H_0^1(\Omega)\cap C(\bar\Omega)$.
Now let $\max_{x\in\Omega}|e^m(x)|=|e^m(x^*)|$ for some $x^*\in\Omega$ (where $x^*$ depends on~$m$).
Also, let $\varsigma^m:={\rm sign}\bigl(e^m(x^*)\bigr)$
and note that $\varsigma^me^m(x^*)=|e^m(x^*)|=\|e^m\|_{L_\infty(\Omega)}$.
Now multiply equation \eqref{semidsicr_em_eq} at $x=x^*$ by $\varsigma^m$ and note that
\begin{subequations}\label{semi_aux}\begin{align}
\varsigma^m(\delta_t^\alpha e^m)\bigr|_{x=x^*}\ge\kappa_{m,m}|e^m(x^*)|-\sum_{j=0}^{m-1}\underbrace{\kappa_{m,j}}_{>0} \|e^j\|_{L_\infty(\Omega)}&= \delta_t^\alpha \|e^m\|_{L_\infty(\Omega)},
\\
\varsigma^m[f(\cdot,t_m,U^m+e^m)-f(\cdot,t_m,U^m)]\bigr|_{x=x^*}\ge -\lambda |e^m(x^*)|&=-\lambda \|e^m\|_{L_\infty(\Omega)},
\end{align}
\end{subequations}
where we used \eqref{delta_def_kappa} and condition A1 on $f$.
Hence, one gets
\beq\label{smooth_e_m}
\varsigma^m\LL e^m(x^*)+(\delta_t^\alpha-\lambda) \|e^m\|_{L_\infty(\Omega)}\le \|r^m\|_{L_\infty(\Omega)}\,.
\eeq
Furthermore, if $e^m$ is sufficiently smooth and $\LL e^m(x^*)$ is defined in the classical sense,
then
$c\ge 0$ implies
$\varsigma^m\LL e^m(x^*)\ge 0$, so \eqref{semi_error} for $p=\infty$ follows.
For less smooth $e^m$, see Remark~\ref{rem_max_princ}.
%
\end{proof}
\smallskip

\begin{remark}[\eqref{semi_error} for $p=\infty$]\label{rem_max_princ}
For less smooth $e^m\in H_0^1(\Omega)\cap C(\bar\Omega)$ in the proof of Theorem~\ref{theo_semidiscr},
split $\Omega$ into disjoint sets $\Omega^+$ and $\Omega^-$ such that $\Omega^\pm:=\{\pm e^m>0\}$.
Next, instead of
\eqref{smooth_e_m},  use similar, but more general relations 
$$
\pm[\LL+\kappa_{m,m}]e(x)\le \sum_{j=0}^{m-1}\kappa_{m,j}\|e^j\|_{L_\infty(\Omega)}+\lambda \|e^m\|_{L_\infty(\Omega)}+\|r^m\|_{L_\infty(\Omega)}
\quad\forall\,x\in\Omega^\pm.
$$
(The above are obtained using \eqref{semi_aux} with $x^*$ replaced by $x\in\Omega^\pm$ and $\varsigma^m$ by $\pm$.)
The desired bound \eqref{semi_error} for $p=\infty$ follows
in view of
\beq\label{max_pr_kappa}
\kappa_{m,m}\|e^m\|_{L_\infty(\Omega^\pm)}\le \sup_{\Omega^\pm}\bigl\{\pm[\LL+\kappa_{m,m}] e^m\bigr\}.
\eeq
The latter is obtained using the maximum principle for functions in $H_0^1(\Omega^\pm)$
\cite[\S8.1]{GTru}.
To be more precise,
note that $\|e^m\|_{L_\infty(\Omega^\pm)}=\sup_{\Omega^\pm}(\pm e^m)$, while the operator $\LL+\kappa_{m,m}$ is linear, so it suffices to get \eqref{max_pr_kappa} only for $\Omega^+$.
Set $M:=\sup_{\Omega^+} \bigl\{[\LL+\kappa_{m,m}]e\bigr\}$. Then by the maximum principle in $\Omega^+$, $M\ge 0$. Consequently, $c\ge 0$
implies that $[\LL+\kappa_{m,m}](M-\kappa_{m,m} e^m)
\ge 0$ in $\Omega^+$, so another application of the maximum principle yields $M-\kappa_{m,m} e^m\ge 0$,
which immediately yields \eqref{max_pr_kappa} for $\Omega^+$. Thus, \eqref{semi_error} is proved.
\end{remark}

\section{Maximum norm error analysis for
finite difference discretizations}\label{sec_FD}

Consider our problem \eqref{problem},\,\eqref{LL_def} in the spatial domain $\Omega=(0,1)^d\subset\R^d$.
Let
 $\bar\Omega_h$ be the tensor product of $d$ uniform meshes $\{ih\}_{i=0}^N$,
 with $\Omega_h:=\bar\Omega_h\backslash\pt\Omega$ denoting the set of interior mesh nodes.
Now, consider a finite difference discretization in the form \eqref{FD_problem},
where $\delta_t^\alpha$ is defined by \eqref{delta_def}. Let the discrete spatial operator $\LL_h$ in \eqref{FD_problem}
be a standard finite difference operator defined,
using the standard orthonormal basis $\{\mathbf{i}_k\}_{k=1}^d$ in $\R^d$
(such that $z=(z_1,\ldots,z_d)=\sum_{k=1}^d z_k\, \mathbf{i}_k$ for any $z\in\R^d$), by
\begin{align*}
&\LL_hV(z):=\\[-0.1cm]
&\!\sum_{k=1}^d h^{-2}\Bigl\{a_k(z+{\textstyle \frac12}h\mathbf{i}_k,t_m)\bigl[V(z)\hspace{-1pt}-\hspace{-1pt}V(z+h\mathbf{i}_k)\bigr]\hspace{-1.3pt}+\hspace{-1pt}a_k(z-{\textstyle \frac12}h\mathbf{i}_k,t_m)\bigl[V(z)\hspace{-1pt}-\hspace{-1pt}V(z-h\mathbf{i}_k)\bigr]\hspace{-1pt}\Bigr\}\\[-0.4cm]
&\qquad{}+\sum_{k=1}^d{\textstyle \frac12}h^{-1}\, b_k(z,t_m)\,\bigl[V(z+h\mathbf{i}_k)-V(z-h\mathbf{i}_k)\bigr] +c(z,t_m)\,V(z)
\qquad\mbox{for}\;\;z\in\Omega_h.
\end{align*}
(Here the terms in the first and second sums respectively discretize $-\pt_{x_k}\!(a_k\,\pt_{x_k}\!u)$ and $b_k\, \pt_{x_k}\!u$ from \eqref{LL_def}.)
The error of this method will be bounded in the nodal maximum norm, denoted
$\|\cdot\|_{L_\infty(\Omega_h)}:=\max_{\Omega_h}|\cdot|$.

We shall assume that $h$ is sufficiently small so that $\LL_h $ satisfies the discrete maximum principle:
\beq\label{d_max_pr}\color{red}
h^{-1}\ge \max_{k=1,\ldots,d}\bigl\{{\textstyle\frac12}\|b_k(\cdot,t)\|_{L_\infty(\Omega)}\,\|a_k(\cdot,t)^{-1}\|_{L_\infty(\Omega)}\bigr\}
\qquad\forall t\in(0,T].
\eeq
Hence, the spatial discrete operator $\LL_h$ is associated with an M-matrix, so
Lemma~\ref{upper_lemma_FD} applies to our discretization.
\smallskip

\begin{theorem}\label{theo_FD}
(i)
Let the temporal mesh satisfy
 \eqref{t_grid_gen}
with $r\ge 1$, and let $\lambda\tau_j^{\alpha}<\{\Gamma(2-\alpha)\}^{-1}$ $\forall\,j\ge1$.
Suppose that $u$ is a unique solution of \eqref{problem},\eqref{LL_def} in $\Omega=(0,1)^d$ with the initial condition $u_0\in L_\infty(\Omega)$ and
 under assumption {\rm A1} on $f$.
Also,
suppose that $\|\partial_t^l u (\cdot, t)\|_{L_\infty(\Omega)}\lesssim 1+t^{\alpha-l}$ for $l = 1,2$ and $t\in(0,T]$,
and
$\|\pt_{x_k}^l u(\cdot,t)\|_{L_\infty(\Omega)}\lesssim 1$ for $l=3,4$, $k=1,\ldots,d$ and $t\in(0,T]$.
Then, under condition
\eqref{d_max_pr} on the above $\LL_h$,
there exists a unique solution $\{U^j\}_{j=0}^M$ of \eqref{FD_problem}, and
\beq\label{FD_error_bound}
\|u(\cdot,t_m)-U^m\|_{L_\infty(\Omega_h)}\lesssim {\mathcal E}^m+t_m^{\alpha}\, h^2
\qquad\forall\,m=1,\ldots,M, 
\eeq
where ${\mathcal E}^m$ is from \eqref{E_cal_m}.

(ii)
If, additionally, $f$ satisfies {\rm A2},
and $\sigma_1\le u_0\le\sigma_2$,
then $\sigma_1\le U^m\le\sigma_2$ $\forall\, m\ge 0$.
\end{theorem}
\smallskip

\begin{proof}
We imitate the proof of Theorem~\ref{theo_semidiscr}.
The existence of a unique solution $\{U^m\}$ follows from Lemma~\ref{upper_lemma_FD}(i), while part (ii)
follows from Lemma~\ref{upper_lemma_FD}(ii),(iii).

It remains to establish~\eqref{FD_error_bound}.
For the error $e^m:= u(\cdot,t_m)-U^m$,  we get a version of \eqref{semidsicr_em_eq} in $\Omega_h$ (instead of $\Omega$) with $\LL$ replaced by $\LL_h$
and $r^m$ replaced by $r^m+r^m_h$,
where $r^m_h:=(\LL_h-\LL)u(\cdot,t_m)$ is the truncation error associated with the spatial discretization.
For the latter, a standard calculation yields $|r^m_h|\lesssim h^2$.

Next, we get the following version of \eqref{semi_error}:
\beq\label{semi_error_FD}
(\delta_t^\alpha-\lambda) \|e^m\|_{L_\infty(\Omega_h)}\le \|r^m+r^m_h\|_{L_\infty(\Omega_h)}
\qquad\forall\,m\ge 1.
\eeq
The proof of the latter closely imitates the proof of \eqref{semi_error} for $p=\infty$, only now
$x^*\in\Omega_h$ is such that
$\max_{x\in\Omega_h}|e^m(x)|=|e^m(x^*)|$, and
a version of \eqref{semi_aux} holds true with $\Omega$ replaced by $\Omega_h$, $\LL$ by $\LL_h$, and $r^m$ by $r^m+r^m_h$.
Finally, in view of \eqref{d_max_pr} combined with $c\ge 0$, one gets
$\varsigma^m\LL_h e^m(x^*)\ge 0$, and hence \eqref{semi_error_FD}.

Let $E^0=E^0_h=0$ and also $(\delta_t^\alpha-\lambda)E^m=\|r^m\|_{L_\infty(\Omega_h)}$
and
$(\delta_t^\alpha-\lambda)E_h^m=\|r_h^m\|_{L_\infty(\Omega_h)}\lesssim h^2$.
Then, applying Corollary~\ref{cor_comparison} to \eqref{semi_error_FD}, one gets $\|e^m\|_{L_\infty(\Omega_h)}\le E^m+E_h^m$.
Also, exactly as in the proof of  Theorem~\ref{lem_simplest_star},  $E^m\lesssim {\mathcal E}^m$.
For $E_h^m$, in view of Theorem~\ref{theo_new_barrier}, the stability property \eqref{main_stab} with $\gamma=-1$ yields
$E_h^m\lesssim t_m^\alpha \,h^2 $. Combining these findings, one gets \eqref{FD_error_bound}.
\end{proof}

\section{Error analysis for finite element discretizations}\label{sec_FE}

Throughout this section, we  restrict our consideration to
the case
$\LL=-\triangle u=-\sum_{k=1}^d\pt^2_{x_k}$
(i.e. $a_k=1$, $b_k=0$ for $k=1,\ldots,d$ and  $c=0$ in \eqref{LL_def}).
Then we
discretize \eqref{problem}, posed in a general bounded  Lipschitz domain  $\Omega\subset\R^d$,
 by applying
 a standard finite element spatial approximation to the temporal semidiscretization~\eqref{semi_semidiscr_method}.
 Let $S_h \subset H_0^1(\Omega)\cap C(\bar\Omega)$ be a Lagrange finite element space of fixed degree $\ell\ge 1$ 
 relative to a
 quasiuniform simplicial triangulation
 $\mathcal T$ of $\Omega$.
 (To simplify the presentation, it will be assumed that the triangulation covers $\Omega$ exactly.)
 Now, for $m=1,\ldots,M$, let $u^m_h \in S_h$ satisfy
 \beq\label{FE_problem}
\begin{array}{l}
\langle \delta_t^\alpha u_h^m,v_h\rangle_h +\langle\nabla u_h^m,\nabla v_h\rangle+ \langle f(\cdot,t_m,u_h^m),v_h\rangle_h=0\qquad\forall v_h\in S_h
\end{array}
\eeq
with $u_h^0= u_0$. 
Here $\langle \cdot, \cdot\rangle$
and $\langle \cdot, \cdot\rangle_h$ respectively denote the exact $L_2(\Omega)$ inner product
and, possibly, its quadrature  approximation.

Our error analysis will employ the standard Ritz projection $\RR_h u(t)\in S_h$ of $u(\cdot,t)$
defined by
\beq\label{Ritz_def}
\langle \nabla\RR_h u,\nabla v_h\rangle=\langle -\triangle u,v_h\rangle_h\qquad\forall v_h\in S_h,
\quad t\in[0,T].
\eeq

\subsection{Lumped-mass linear finite elements:  error analysis in the $L_\infty(\Omega)$ norm}\label{ssec_lumped}
First, we consider
lumped-mass linear finite-element discretizations,
i.e.
$\ell=1$ and
$\langle \cdot,\cdot\rangle_h$ in \eqref{FE_problem} is defined using the quadrature rule $Q_T[v]:=\int_T v^I$, where $v^I$ is the standard linear Lagrange interpolant.

Let $\mathcal N$ denote the set of interior mesh nodes,  with the corresponding  piecewise-linear basis hat functions $\{\phi_z\}_{z\in\mathcal N}$.
Then,
 using $v_h=\phi_z$ in \eqref{FE_problem}, our discretization can be represented in the form of the discrete problem \eqref{FD_problem} for the nodal values of the computed solution $U^m(z):=u^m_h(z)$, with
\beq\label{LLh_fem}
\LL_h U^m(z):=\frac{\langle \nabla u^m_h,\nabla \phi_z\rangle}{\langle 1, \phi_z\rangle_h}
\qquad\forall\,z\in{\mathcal N}.
\eeq
We shall additionally assume that the spatial triangulation is such that $\LL_h$ is associated with an M-matrix
(sufficient conditions for this are discussed in Remark~\ref{rem_Del_triang}).
Hence, Lemma~\ref{upper_lemma_FD} applies to our finite element discretization.

Our main result for this discretization is the following.
\smallskip

\begin{theorem}[Lumped-mass linear elements]\label{theo_lumped_FE}
(i)
Let the temporal mesh satisfy
 \eqref{t_grid_gen}
with $r\ge 1$, and let $\lambda\tau_j^{\alpha}<\{\Gamma(2-\alpha)\}^{-1}$ $\forall\,j\ge1$.
Suppose that $u$ is a unique solution of \eqref{problem},\eqref{LL_def}  with the initial condition $u_0\in L_\infty(\Omega)$ and
 under assumption {\rm A1} on $f$.
Also,
suppose that$\|\partial_t^l u (\cdot, t)\|_{L_\infty(\Omega)}\lesssim 1+t^{\alpha-l}$ for $l = 1,2$
and $t\in(0,T]$.
%
Then, if the operator $\LL_h$ from \eqref{LLh_fem} is associated with an M-matrix,
there exists a unique solution $\{u_h^j\}_{j=0}^M$ of \eqref{FE_problem}, and, for $m=1,\ldots,M$,
\beq
\|u(\cdot,t_m)-u_h^m\|_{L_\infty(\Omega)}\lesssim
{\mathcal E}^m
 \label{FElin_error_bound}
+\max_{t\in\{0,t_m\}}\|\rho(\cdot, t)\|_{L_\infty(\Omega)}+\int_0^{t_m}\!\|\pt_t \rho(\cdot, t)\|_{L_\infty(\Omega)}\,dt,
\eeq
where ${\mathcal E}^m$ is defined in \eqref{E_cal_m}, and
 $\rho(\cdot, t):=\RR_h u(t)-u(\cdot, t)$
 is the error of the Ritz projection~\eqref{Ritz_def}.

(ii)
If, additionally, $f$ satisfies {\rm A2},
and $\sigma_1\le u_0\le\sigma_2$,
then $\sigma_1\le u_h^m\le\sigma_2$ $\forall\, m\ge 0$.
\end{theorem}
\smallskip

\begin{proof}
We imitate the proofs of Theorems~\ref{theo_semidiscr} and~\ref{theo_FD}.
First, since
our discretization can be represented in the form of the discrete problem \eqref{FD_problem} for the nodal values of the computed solution $U^m(z)=u^m_h(z)$,
the existence of a unique solution $\{U^m\}$ follows from Lemma~\ref{upper_lemma_FD}(i), while part (ii)
follows from Lemma~\ref{upper_lemma_FD}(ii),(iii).

It remains to establish~\eqref{FElin_error_bound}.
Note that
$u(\cdot,t_m)-u_h^m=[\RR_h u(\cdot,t_m)-u_h^m]-\rho(\cdot, t_m)$, where $\RR_h u(\cdot,t_m)-u_h^m\in S_h$.
Hence, it suffices to prove the desired bound for the nodal values of the latter, which will be denoted by
 $e^m:= \RR_h u(\cdot,t_m)-U^m$  $\forall\,z\in\mathcal N$.

In view of \eqref{LLh_fem}, one has
$\LL_h \RR_h u(z,t_m)=-\triangle u(z,t_m)$. Or, equivalently,
using \eqref{problem} and the truncation error $r^m=\delta_{t}^{\alpha} u(\cdot,t_m)-D_t^\alpha u(\cdot,t_m)$,
one can rewrite it as
$$
\delta^\alpha_tu(z,t_m)+\LL_h \RR_h u(z,t_m)+f(z,t_m,u(z,t_m))=r^m\qquad \forall\,z\in\mathcal N,\;\forall\,m\ge 1.
$$
Subtracting the nodal representation \eqref{FD_problem} of  our discretization,
one gets
$e^0=\rho^0$ and
\beq\label{error_eq_simple_lumptedFE}
\delta_t^\alpha[ e^m-\rho^m] +\LL_h e^m+[f(\cdot, t_m,u(\cdot,t_m))-f(\cdot, t_m,U^m)]=
r^m
\quad \forall\,z\in\mathcal N,\;\forall\,m\ge 1
\eeq
(which is a version of \eqref{error_eq_simple}),
where we used the notation $\rho^m:=\rho(\cdot,t_m)$ at any $z\in\mathcal N$.
Next, using the constant $\lambda\ge 0 $ from assumption A1 on $f$, set
$$
p^m
:=\lambda +\left\{\!\!
\begin{array}{cl}
\frac{f(\cdot, t_m,u(\cdot,t_m))-f(\cdot, t_m,U^m)}{u(\cdot,t_m)-U^m},&
\mbox{if~}u(\cdot,t_m)\neq U^m,\\
0,&\mbox{otherwise,}
\end{array}
\right.
\quad\forall\,z\in\mathcal N,\;\forall\,m\ge 1.
$$
Then, in view of A1, $p^m\ge 0$.
Also, $f(\cdot, t_m,u(\cdot,t_m))-f(\cdot, t_m,U^m)=(p^m-\lambda)[e^m-\rho^m]$,
so \eqref{error_eq_simple_lumptedFE} can be rewritten as
\beq\label{error_eq_simple_lumptedFE_new}
(\delta_t^\alpha+\LL_h +p^m-\lambda) e^m=
r^m+(p^m-\lambda)\rho^m+\delta_t^\alpha \rho^m
\qquad \forall\,z\in\mathcal N,\;\forall\,m\ge 1.
\eeq
This is a linear version of \eqref{FD_problem}, so, on the one hand,
in view of Lemma~\ref{upper_lemma_FD}(ii),(iii),
we can construct upper and lower solutions to estimate $e^m$.
On the other hand, we can separately estimate the components of the error that correspond to the three terms in the right-hand side of \eqref{error_eq_simple_lumptedFE_new}.

First, suppose that the right-hand side of \eqref{error_eq_simple_lumptedFE_new} equals $r^m$ and $e^0=0$.
Then for $E^m$ such that $E^0=0$ and also $(\delta_t^\alpha-\lambda)E^m=\|r^m\|_{L_\infty(\Omega)}$, exactly as in the proof of  Theorem~\ref{lem_simplest_star}, one gets $E^m\lesssim {\mathcal E}^m$.
Also, by \eqref{LLh_fem}, $(\LL_h +p^m) E^m=p^m E^m\ge 0$.
Hence, the pair
$\pm E^m$ gives discrete upper and lower solutions for \eqref{error_eq_simple_lumptedFE_new} in this case. So
$|e^m|\le E^m\lesssim {\mathcal E}^m$, and the desired bound of type \eqref{FElin_error_bound} on $\|e^m\|_{L_\infty(\Omega)}$ follows.

Next, suppose that $e^0=\rho^0$ and the right-hand side of \eqref{error_eq_simple_lumptedFE_new} equals $(p^m-\lambda)\rho^m$
(where no upper bound on $p^m$ is available).
Let $B^0=0$ and $(\delta_t^\alpha-\lambda)B^m=1$, so, in view of Theorem~\ref{theo_new_barrier}, $0\le B^m\lesssim t_m^\alpha$.
Next, note that $(\delta_t^\alpha-\lambda)[2\lambda B^m+1]=\lambda$,
while $(\LL_h+p^m)[2\lambda B^m+1]=p^m [2\lambda B^m+1]\ge p^m$.
Consequently, the pair of functions
$\pm[2\lambda B^m+1]\sup_{[0,t_M]}\|\rho\|_{L_\infty(\Omega)}$ gives discrete upper and lower solutions
for \eqref{error_eq_simple_lumptedFE_new} in this case.
Hence, $ |e^m|\le [2\lambda B^m+1]\sup_{[0,t_M]}\|\rho\|_{L_\infty(\Omega)}$, so one immediately gets
$\|e^M\|_{L_\infty(\Omega)}\lesssim \sup_{[0,t_M]}\|\rho\|_{L_\infty(\Omega)}$.
As a similar argument applies for any $M\ge1$, we deduce the desired bound of type \eqref{FElin_error_bound} on $\|e^m\|_{L_\infty(\Omega)}$.

In a similar manner, consider \eqref{error_eq_simple_lumptedFE_new} with the right-hand side equal to $\delta_t^\alpha \rho^m$ and $e^0=0$.
Let
$\bar\rho^m:=\int_0^{t_m}\|\pt_t \rho(\cdot, s)\|_{L_\infty(\Omega)}\,ds$, for which,
in view of \eqref{delta_def},
one gets
$|\delta_t\rho^m|\le \delta_t\bar\rho^m$, and so $|\delta^\alpha_t\rho^m|\le \delta^\alpha_t\bar\rho^m$.
Consequently, the pair of functions
$\pm[\bar\rho^m+\lambda\bar\rho^M B^m]$ gives discrete upper and lower solutions
for \eqref{error_eq_simple_lumptedFE_new} in this case.
Hence, $\|e^M\|_{L_\infty(\Omega)}\lesssim \bar\rho^M$.
Applying a similar argument  for any $M\ge1$, we again deduce the desired bound of type \eqref{FElin_error_bound} on $\|e^m\|_{L_\infty(\Omega)}$.
\end{proof}
\smallskip

\begin{remark}[$\LL_h$ associated with an M-matrix]\label{rem_Del_triang}
The operator $\LL_h$ from \eqref{LLh_fem} is associated
with a normalized stiffness matrix for $-\triangle$. The latter is
 an M-matrix under the following  conditions
 on the triangulation.
For $\Omega\subset\R^2$, let $\mathcal T$ be a Delaunay triangulation, i.e.,
the sum of the angles opposite to any interior
edge is less than or equal to $\pi$.
In the case $\Omega\subset\R^3$,
it is sufficient, but not necessary, for the triangulation to be non-obtuse (i.e. with no interior angle in any mesh element exceeding $\frac{\pi}2$).
For weaker necessary and sufficient conditions, we refer the reader to \cite[Lemma~2.1]{xu_zik}.
\end{remark}
\smallskip

\begin{remark}[Ritz projection]\label{rem_ritz}
The error bound \eqref{FElin_error_bound} involves $\rho$, the error of the Ritz projection.
For the latter, assuming that the spatial domain $\Omega$ is polygonal,
convex polyhedral or smooth,
for the considered lumped-mass discretization,
one has \cite[(5.6)]{NK_MC_L1}\vspace{-0pt}
$$
\|\pt_t^l \rho(\cdot, t)\|_{L_\infty(\Omega)}\lesssim h^{2-q}|\ln h|\Bigl\{\|\pt_t^l u(\cdot, t)\|_{W^{2-q}_\infty(\Omega)}
+
\|\pt_t^l \LL u(\cdot, t)\|_{W^{2-q}_{d/2}
(\Omega)}\Bigr\},
$$
where $l=0,1$, $q=0,1$ and $t\in(0,T]$.
Thus, under certain realistic assumptions on $u$ (see, e.g., \cite[Corollary~5.7 and Remark~5.8]{NK_MC_L1}),
the error bound~\eqref{FElin_error_bound} yields
$\|u(\cdot,t_m)-u_h^m\|_{L_\infty(\Omega)}\lesssim{\mathcal E}^m+h^2|\ln h|$.
\end{remark}

\subsection{Finite elements without quadrature:  error analysis in the $L_2(\Omega)$ norm}\label{ssec_L1_full_FEM}
Next, consider
finite elements of fixed degree $\ell\ge 1$ without quadrature, i.e. with $\langle \cdot, \cdot\rangle_h=\langle \cdot, \cdot\rangle$
in \eqref{FE_problem}.
We shall need an additional assumption on $f$.
\medskip

\begin{itemize}
\item[{\bf A1$\mathbf{}^*$\!\!.}]
Let $f$ satisfy the one-sided Lipschitz condition
$$
|f(x,t,s_1)-f(x,t,s_2)|\le \bar\lambda|s_1-s_2|\qquad \forall s_1,\, s_2\in\R,\;\; x\in\Omega,\;\;t>0
$$
with some constant $\bar\lambda\ge0$. (Clearly, $\bar\lambda\ge\lambda$ for $\lambda$ from A1.)
\end{itemize}
\medskip

\begin{theorem}
\label{theo_FE_L2norm}
Let the temporal mesh satisfy
 \eqref{t_grid_gen}
with $r\ge 1$, and let $\lambda\tau_j^{\alpha}<\{\Gamma(2-\alpha)\}^{-1}$ $\forall\,j\ge1$.
Suppose that $u$ is a unique solution of \eqref{problem},\eqref{LL_def}  with the initial condition $u_0\in L_\infty(\Omega)$ and
 under assumptions {\rm A1} and {\rm A1${}^*$} on $f$.
Also,
suppose that $\|\partial_t^l u (\cdot, t)\|_{L_2(\Omega)}\lesssim 1+t^{\alpha-l}$ for $l = 1,2$
and $t\in(0,T]$.
%
Then, under the condition $\langle \cdot, \cdot\rangle_h=\langle \cdot, \cdot\rangle$,
there exists a unique solution $\{u_h^j\}_{j=0}^M$ of \eqref{FE_problem}, and, for $m=1,\ldots,M$,
\beq
\|u(\cdot,t_m)-u_h^m\|_{L_2(\Omega)}\lesssim
{\mathcal E}^m
 \label{FE_error_bound_L2}
+\max_{t\in[0,t_m]}\|\rho(\cdot, t)\|_{L_2(\Omega)}+\int_0^{t_m}\!\|\pt_t \rho(\cdot, t)\|_{L_2(\Omega)}\,dt,
\eeq
where ${\mathcal E}^m$ is defined in \eqref{E_cal_m}, and
 $\rho(\cdot, t):=\RR_h u(t)-u(\cdot, t)$
 is the error of the Ritz projection~\eqref{Ritz_def}.
\end{theorem}
\smallskip

\begin{proof}
The existence of a unique solution $u_h^m$ is established
noting that, in view of A1 and the upper bound on $\lambda\tau_j^{\alpha}$, at each time level $t_m$ we have a finite element discretization of type \eqref{FE_problem} for the monotone elliptic equation
\eqref{semi_elliptic} (as discussed in the proof of Lemma~\ref{upper_lemma}(i)).
Hence, the latter finite element discretization  is equivalent to the minimization of a uniformly convex and continuously differentiable functional on a finite-dimensional space, so the existence of a unique computed solution
follows (see, e.g., \cite[\S4.3.9]{orgega_r_book}).

It remains to obtain the error bound \eqref{FE_error_bound_L2}, for which we shall partially imitate the proofs of
 Theorems~\ref{theo_semidiscr} and~\ref{theo_lumped_FE}.
Let $e_h^m:=\RR_h u(t_m)-u_h^m\in S_h$ and $\rho^m:=\rho(\cdot, t_m)$.
Then $u(\cdot,t_m)-u_h^m=e_h^m-\rho^m$, so
it suffices to prove the desired bounds for $e_h^m$.
Now, a standard calculation using \eqref{FE_problem} (in which $\langle \cdot, \cdot\rangle_h=\langle \cdot, \cdot\rangle$) and \eqref{problem} yields
\beq\label{FE_err_prob_L2}
\langle \delta_t^\alpha e_h^m, v_h\rangle +\langle\nabla e_h^m,\nabla v_h\rangle
+\langle f(\cdot,t_m, u(\cdot,t_m))-f(\cdot,t_m u_h^m), v_h\rangle=\langle \delta_t^\alpha \rho^m+r^m, v_h\rangle
\eeq
$\forall\, v_h\in S_h$.
Here we again use the truncation error $r^m=\delta_t^\alpha u(\cdot,t_m)-D_t^\alpha u(\cdot,t_m)$,
for which we again have \eqref{trunc_er_semidis} with $p=2$.
Next, note that $u(\cdot,t_m)=u_h^m+e_h^m-\rho^m$.
So, setting $v_h:=e_h^m$ and recalling A1${}^*$, we arrive at
$$
\langle \delta_t^\alpha e_h^m, e_h^m\rangle
+\langle f(\cdot,t_m, u_h^m+e_h^m)-f(\cdot,t_m u_h^m), e_h^m\rangle\le \langle r^m+\delta_t^\alpha \rho^m,e_h^m\rangle+\bar\lambda\langle|\rho^m|,|e_h^m|\rangle.
$$
The left-hand side here is estimated using a version of \eqref{semi_aux_p2} (with $e^m$ replaced by $e^m_h$ and $U^m$ replaced by $u_h^m$). Hence, we get the following version of
\eqref{semi_error}:
\beq\label{eq_e_m_L2}
(\delta_t^\alpha-\lambda) \|e^m\|_{L_2(\Omega)}\le \|r^m\|_{L_2(\Omega)}
+\|\delta_t^\alpha \rho^m\|_{L_2(\Omega)}+\bar\lambda\|\rho^m\|_{L_2(\Omega)}
\qquad\forall\,m\ge 1,
\eeq
subject to $e_h^0=\rho^0$.

Let $E^0=B^0=0$, and also $(\delta_t^\alpha-\lambda)E^m=\|r^m\|_{L_2(\Omega)}$ and $(\delta_t^\alpha-\lambda)B^m=1$.
Then, exactly as in the proof of  Theorem~\ref{lem_simplest_star}, one gets $E^m\lesssim {\mathcal E}^m$.
Also, in view of Theorem~\ref{theo_new_barrier}, $0\le B^m\lesssim t_m^\alpha$.
Additionally, consider
$\bar\rho^m:=\int_0^{t_m}\|\pt_t \rho(\cdot, s)\|_{L_2(\Omega)}\,ds$, for which,
in view of \eqref{delta_def},
one gets
$\|\delta_t\rho^m\|_{L_2(\Omega)}\le \delta_t\bar\rho^m$, and so $\|\delta^\alpha_t\rho^m\|_{L_2(\Omega)}\le \delta^\alpha_t\bar\rho^m$.
Consequently, $(\delta_t^\alpha-\lambda)\bar\rho^m\ge \|\delta^\alpha_t\rho^m\|_{L_2(\Omega)}-\lambda \bar\rho^M$.
Combining these findings, one concludes that the function
$$
\|\rho^0\|_{L_2(\Omega)}+E^m+\bar\rho^m+\Bigl(\lambda \|\rho^0\|_{L_2(\Omega)}+\lambda \bar\rho^M +\bar\lambda\max_{j=0,\ldots,M}\|\rho^j\|_{L_2(\Omega)}\Bigr)B^m
$$
is an upper solution for problem \eqref{eq_e_m_L2}.
Hence, in view of Corollary~\ref{cor_comparison}, one gets the desired bound~\eqref{FE_error_bound_L2} for $m=M$.
Applying a similar argument  for any $M\ge1$, we again deduce the desired bound $\forall m\ge1$.
\end{proof}
\smallskip

\begin{remark}[Ritz projection]\label{rem_ritz2}
The error bound \eqref{FE_error_bound_L2} involves $\rho$, the error of the Ritz projection.
For the latter, assuming that the spatial domain $\Omega$ is smooth or convex (or, more generally, such that $\|v\|_{W^2_2(\Omega)}\lesssim \|\LL v\|_{L_2(\Omega)}$ for any sufficiently smooth $v$),
one has
$$
\|\pt_t^l \rho(\cdot, t)\|_{L_2(\Omega)}\lesssim h\inf_{v_h\in S_h}\|\pt_t^l u(\cdot, t)-v_h\|_{W^1_2(\Omega)}
\quad\;\;\mbox{for}\;\;
l=0,1,\;t\in(0,T].
$$
For $l=0$, see, e.g., \cite[Theorem~5.7.6]{BrenScott}.
A similar result for $l=1$ follows as $\pt_t\rho(\cdot, t)=\RR_h \dot u(t)-\dot u(\cdot, t)$, where $\dot u:=\pt_t u$.
Thus, under certain realistic assumptions on $u$ (see, e.g., \cite[Corollary~5.3 and Remark~5.4]{NK_MC_L1}),
\eqref{FE_error_bound_L2} yields
$\|u(\cdot,t_m)-u_h^m\|_{L_2(\Omega)}\lesssim{\mathcal E}^m+h^{\ell+1}$.
\end{remark}

\begin{remark}[More general $\LL$]\color{red}
  Theorem~\ref{theo_FE_L2norm} can be immediately extended to the case of more general $\LL=\LL(t)$
 associated with a coercive bilinear form.
    The only modification required in the proof is
    to replace $\langle\nabla e_h^m,\nabla v_h\rangle$
     in \eqref{FE_err_prob_L2} by $\langle \LL (t_m)e_h^m,v_h\rangle$. As $\langle \LL (t_m)e_h^m,e_h^m\rangle\ge 0$, we again get \eqref{eq_e_m_L2}, so the remainder of the proof works without any further changes.
    Note that the estimation of the error of the Ritz projection (such as discussed in Remark~\ref{rem_ritz2})
     will be more complicated in this case.
\end{remark}

\section{Generalizations}\label{sub_gen}

\subsection{A2 satisfied, but not A1}\label{ssub_gen_A1}
Suppose that
 $f$ in \eqref{problem} satisfies {\rm A2}, but not A1 (as, e.g., in the Fisher equation with $f=u^2-u$), and the initial condition is such that
 $\sigma_1\le u_0\le \sigma_2$.
 Also,
 let $f$ be continuous in $s$ and  satisfy $f(\cdot,t,s)\in L_\infty(\Omega)$ for all $t>0$ and $s\in[\sigma_1,\sigma_2]$.

Then one can replace $f$ with a standard modification $\widetilde f=\widetilde f(\cdot,t,s)$ defined by
$\widetilde f:=f$ for $s\in[\sigma_1,\sigma_2]$, and $\widetilde f=f(\cdot,t,\sigma_1)$ for $s\le\sigma_1$, and
$\widetilde f=f(\cdot,t,\sigma_2)$ for $s\ge\sigma_2$.
Clearly $\widetilde f$ satisfies both A1 and A2, as well as A1${}^*$, so all our results on existence, uniqueness and convergence properties of the discrete solutions will apply. Furthermore, with the exception of Theorem~\ref{theo_FE_L2norm},
the computed solutions will lie between $\sigma_1$ and $\sigma_2$; hence they will also be (not necessarily unique) solutions of the corresponding discrete problems with the original $f$.
Note also that the nonlinear discrete problems with $\widetilde f$ may be computationally more stable.

{\color{red}
\subsection{Nonhomogeneous Dirichlet boundary condition}
Suppose that
$u=\varphi$ on $\pt\Omega\times(0,T]$ in \eqref{problem}, where $\varphi(\cdot,t)\in H^1(\Omega)\cap L_\infty(\Omega)$ $\forall\,t\in(0,T]$,
while $\sigma_1\le \varphi\le\sigma_2$ on $\pt\Omega\times(0,T]$ in A2.
Then, with the obvious modifications
$ U^m=\varphi(\cdot,t_m)$  on $\pt\Omega$ in \eqref{semi_semidiscr_method} and
 $\bar U^m\ge\varphi(\cdot, t_m)$ on $\pt\Omega$
  in \eqref{upper_semidiscr_method}, and a similar change in \eqref{FD_problem},
all results of \S\ref{sec_upper} remain valid.
In particular, in the proof of Lemma~\ref{upper_lemma}(i),
the existence  of a unique solution of \eqref{semi_elliptic} such that $U^m-\varphi(\cdot,t_m)\in H^1_0(\Omega)\cap L_\infty(\Omega)$
can be shown imitating the argument used in the final paragraph of \cite[\S2.1]{DK16}.
Furthermore, all error bounds of \S\ref{ssec_L1_semi} and \S\ref{sec_FD} remain valid for this case.
Similarly, the results of \S\ref{sec_FE} 
can also be extended for this case
with the obvious changes
$u_h^m-\varphi_h(\cdot,t_m)\in S_h$ in \eqref{FE_problem} and ${\mathcal R}_hu-\varphi_h\in S_h$ in \eqref{Ritz_def},
where $\varphi_h$ is a standard Lagrange interpolant of $\varphi$,
 except the bounds on the Ritz projection
in Remarks~\ref{rem_ritz} and~\ref{rem_ritz2} should now take into account the error induced by $\varphi-\varphi_h$.
}

\subsection{Periodic boundary conditions}
As many of our arguments rely on the discrete maximum principle for the spatial operator $\LL_h$, they can easily be extended to other types of boundary conditions.
In particular, the results of \S\ref{sec_FD} for finite difference discretizations in $\Omega=(0,1)^d$,
including Theorem~\ref{theo_FD}, apply
to the case of periodic boundary conditions (with standard modifications in~\eqref{FD_problem} to reflect such boundary conditions). Note that a version  of Lemma~\ref{upper_lemma_FD} from \S\ref{ssec2_2} also holds true for this case assuming that the strict version of
$\lambda\tau_j^{\alpha}\le \{\Gamma(2-\alpha)\}^{-1}$ is satisfied.

\subsection{Neumann/Robin and mixed  boundary conditions}

Suppose that on  a subset $\pt\Omega_R$ of the boundary $\pt\Omega$, the Dirichlet boundary condition in \eqref{problem} is replaced by the homogeneous Neumann/Robin boundary condition of the form
\beq\label{BC_Robin}
\frac{\partial u}{\partial n}+\mu u=0\qquad\mbox{on~~}\pt\Omega_R\subseteq\pt\Omega,
\qquad\mbox{where~~}\mu(x,t)\ge0.
\eeq
Then Lemma~\ref{upper_lemma_FD} from \S\ref{ssec2_2} remains true
provided that $\pt\Omega$ in~\eqref{FD_problem} is replaced by $\pt\Omega\backslash\pt\Omega_R$, so $\Omega_h$ includes the nodes on $\pt\Omega_R$, and also
the strict version of
$\lambda\tau_j^{\alpha}\le \{\Gamma(2-\alpha)\}^{-1}$ is satisfied.
(In fact, the latter is required only if $\pt\Omega_R=\pt\Omega$ and $\mu=0$ on $\pt\Omega$.)
Now, consider the treatment of \eqref{BC_Robin} in finite difference and finite element approximations separately.

\subsubsection{Finite difference discretizations}

The material of \S\ref{sec_FD} can be also extended for \eqref{BC_Robin}.
Using the standard finite difference discretization of the Robin boundary conditions (see, e.g.,
\cite[\S{}VII.1.9]{Sam}), we modify the definition of $\LL_h V(z)$ for $z\in\pt\Omega_R$ as follows.
Whenever $z\in\pt\Omega_R$ and $z\pm h\mathbf{i}_k\not\in\Omega$, we replace
$V(z\pm h\mathbf{i}_k)$ in $\LL_h V(z)$ by $V(z\mp h\mathbf{i}_k)+2h\,\mu(z,t_m)\,V(z)$.

The same condition~\eqref{d_max_pr} ensures that $\LL_h$ satisfies the discrete maximum principle also in this case.
However, we need to modify the proof of Theorem~\ref{theo_FD}, as the truncation error associated with the spatial discretization
$r^m_h=(\LL_h-\LL)u(\cdot,t_m)$ is only $O(h)$ on $\pt\Omega_R$
(while $|r^m_h|\lesssim h^2$ on $\Omega_h\backslash\pt\Omega_R$).

\smallskip
\renewcommand{\thetheoremaa}{\ref{theo_FD}${}^*$}
\begin{theoremaa}
\label{theo_FD_Robin}
Let the coefficients $\{a_k\}$ in \eqref{LL_def} be positive constants, and
$\pt\Omega_R\subseteq \pt\Omega$.
Then Theorem~\ref{theo_FD} holds true for the above finite difference discretization
with $t_m^\alpha h^2$ in the right-hand side of the error bound \eqref{FD_error_bound}
replaced by~$h^2$.%
\end{theoremaa}%
\smallskip

\begin{proof}
Imitating the proof of Theorem~\ref{theo_FD}, we again get the following version of
\eqref{semidsicr_em_eq} in $\Omega_h$ (only now $\Omega_h$ includes the nodes on $\pt\Omega_R$):
\beq\label{FD_er_Robin}
\delta_t^\alpha e^m +\LL_h e^m+[f(\cdot, t_m,U^m+e^m)-f(\cdot, t_m,U^m)]=r^m+r^m_h\qquad \forall\,m\ge 1.
\eeq
Next, similarly to obtaining \eqref{error_eq_simple_lumptedFE_new} in the proof of Theorem~\ref{theo_lumped_FE}, introduce $p^m\ge 0$ such that the above is
rewritten in the form
\beq\label{FD_er_Robin1}
(\delta_t^\alpha+\LL_h +p^m-\lambda) e^m=
r^m+r^m_h\qquad \forall\,m\ge 1.
\eeq
Set $r^m_R:=0$ in $\Omega$ and $r^m_R:=r^m_h=O(h)$ on $\pt\Omega_R$.
As the above is a linear version of \eqref{FD_er_Robin}, we can separately estimate the components of the error that correspond to  $r^m_R$ and $r^m+(r^m_h-r^m_R)$. For the latter, exactly as in the proof of Theorem~\ref{theo_FD}, we get
a version of \eqref{semi_error_FD} with $r^m_h$ replaced by $(r^m_h-r^m_R)=O(h^2)$, so the desired error bound of type \eqref{FD_error_bound} for this component of the error follows.

The remaining component of the error  satisfies \eqref{FD_er_Robin1} with the right-hand side $r^m_R$, and, in view of Lemma~\ref{upper_lemma_FD}(ii),(iii), can be estimated using upper and lower solutions.
To simplify the presentation, we shall assume that $\pt\Omega_R\subset\{x_1=1\}$
{\color{red}and $b_1$ is constant}
(as the other cases are similar).
Let $B^0=0$ and $(\delta_t^\alpha-\lambda)B^m=1$, so, in view of Theorem~\ref{theo_new_barrier}, $0\le B^m\lesssim t_m^\alpha$.
A calculation shows that $\LL_h x_1\ge \color{blue}-|b_1|$ in $\Omega$ (where we exploit that the coefficient $a_1$ is constant),
while $\LL_h x_1\ge 2 h^{-1} a_1$ on $\pt\Omega_R$.
Noting that $(\delta_t^\alpha+\LL_h +p^m-\lambda)x_1\ge \LL_h x_1-\lambda $, one can check that
the pair of discrete functions
$$
\pm h(2a_1)^{-1}[x_1+{\color{red}(\lambda+|b_1|)} B^m]\max_{m=1,\ldots,M}\|r^m_R\|_{L_\infty(\pt\Omega_R)}
$$
gives an upper and a lower solutions for the component of $e^m$ that we are estimating.
As $\|r^m_R\|_{L_\infty(\pt\Omega_R)}\lesssim h$, we conclude that this component of the error
is $\lesssim h^2$.
\end{proof}

\subsubsection{Lumped-mass linear finite elements}\label{subsub_FE}
Next, consider an extension of the material of \S\ref{ssec_lumped} for \eqref{BC_Robin}.
To simplify the presentation, let $\pt\Omega_R\neq\pt\Omega$ or $\mu>0$ (to ensure that the Ritz projection is well-defined).
In this case, with an obvious modification of $S_h$, the standard lumped-mass discretization \eqref{FE_problem} will include an additional term
$\int_{\pt\Omega_R}\bigl(\mu(\cdot, t_m)\,u_h^m\, v_h\bigr){}^I$ in the left-hand side.
A similar modification applies  to the definition of the Ritz projection \eqref{Ritz_def},
in which the left-hand side now includes an additional term $\int_{\pt\Omega_R}\bigl(\mu(\cdot, t_m)\, v_h\,\RR_h u\bigr){}^I$.
Finally, in the definition of $\LL_h$ in \eqref{LLh_fem}, the term
$\langle \nabla u^m_h,\nabla \phi_z\rangle$ is now replaced by $\langle \nabla u^m_h,\nabla \phi_z\rangle+\int_{\pt\Omega_R}\bigl(\mu(\cdot, t_m)u_h^m \phi_z\bigr){}^I$, while $\mathcal N$
denotes the set of nodes in $\Omega\cup\pt\Omega_R$.
With these modifications, an inspection of the proof of Theorem~\ref{theo_lumped_FE} shows that this theorem remains true.

\subsubsection{Finite elements without quadrature}
Finally, we proceed to an extension of \S\ref{ssec_L1_full_FEM}.
The treatment of the boundary condition \eqref{BC_Robin} remains as in \S\ref{subsub_FE}, only
all approximate integrals of type $\int_{\pt\Omega_R}(\cdots)^I$ are now replaced by their exact versions
$\int_{\pt\Omega_R}(\cdots)$.
Then an inspection of the proof of Theorem~\ref{theo_FE_L2norm} shows that in
\eqref{FE_err_prob_L2} we need to add $\int_{\pt\Omega_R}\mu(\cdot,t_m)\,e_h^m\, v_h$ to
$\langle\nabla e_h^m,\nabla v_h\rangle$, and afterwards, when we set $v_h:=e_h^m$ in \eqref{FE_err_prob_L2},
we now exploit the positivity of $\langle\nabla e_h^m,\nabla e_h^m\rangle$ and $\int_{\pt\Omega_R}\mu(\cdot,t_m)(e_h^m)^2$.
Thus, we conclude that Theorem~\ref{theo_FE_L2norm} remains valid for the considered finite element discretization.

\section{Numerical results}\label{sec_Num}


\begin{table}[t!]
\begin{center}
\caption{%
{
Maximum  nodal errors at $t=1$ (odd rows) and
computational rates $q$ in $M^{-q}$ or $N^{-q}$ (even rows) on the graded mesh with
$r=1$,
$r=(2-\alpha)/.9$ and $r=(2-\alpha)/\alpha$}
}
\label{t_positive_time}
\vspace{-0.1cm}
\tabcolsep=4pt
{\small
\begin{tabular}{lrrrrrrrr}
\hline
\strut\rule{0pt}{9pt}
&\multicolumn{4}{l}{errors and convergence rates in time}&
\multicolumn{4}{l}{~~~errors and convergence rates in space}\\[0pt]
\strut\rule{0pt}{9pt}
&\multicolumn{4}{l}{$N = 2M$}&
\multicolumn{4}{l}{~~~$M = N^2
$}\\[2pt]
\hline
\strut\rule{0pt}{9pt}
&$M=2^5$& $M=2^6$& $M=2^7$& $M=2^8$&{}~~~$N=2^3$&~$N=2^4$&~$N=2^5$&~$N=2^6$\\
\hline

\strut\rule{0pt}{11pt}&\multicolumn{8}{c}{$r=1$}\\[-1pt]\cline{2-9}\rule{0pt}{9pt}

$\alpha=0.3\;\;$
&1.88e-3	&8.98e-4	&4.37e-4	&2.15e-4	
    &1.23e-2	&2.99e-3	&7.49e-4	&1.87e-4\\
&1.07	&1.04	&1.02&  &2.05	&2.00	&2.00\\[3pt]	
$\alpha=0.5$
&7.41e-4	&3.35e-4	&1.58e-4	&7.65e-5
    &8.09e-3	&2.07e-3	&5.13e-4	&1.28e-4\\
&1.15	&1.08	&1.05&  &1.97	&2.01	&2.00\\[3pt]
$\alpha=0.7$ 	  	
&1.06e-3	&4.83e-4	&2.27e-4	&1.08e-4
    &5.87e-3	&1.48e-3	&3.67e-4	&9.14e-5\\
&1.13	&1.09	&1.06&&1.98	&2.02	&2.01\\
       	
\hline

\strut\rule{0pt}{11pt}&\multicolumn{8}{c}{$r=\frac{2-\alpha}{.9}$}\\[1.5pt]\cline{2-9}\rule{0pt}{9pt}

$\alpha=0.3\;\;$
&5.87e-4	&1.79e-4	&5.49e-5	&1.69e-5	
    &1.15e-2	&2.81e-3	&7.04e-4	&1.75e-4\\
&1.71	&1.71	&1.70&&2.04	&2.00	&2.00\\[3pt]	
$\alpha=0.5$
&3.30e-4	&1.09e-4	&3.70e-5	&1.29e-5
    &7.88e-3	&2.01e-3	&4.98e-4	&1.24e-4\\
&1.60	&1.56	&1.53&&1.97	&2.01	&2.00\\[3pt]
$\alpha=0.7$ 	  	
&7.14e-4	&2.83e-4	&1.15e-4	&4.75e-5
    &5.66e-3	&1.42e-3	&3.49e-4	&8.66e-5\\
&1.33	&1.30	&1.28&&1.99	&2.02	&2.01\\

\hline

\strut\rule{0pt}{11pt}&\multicolumn{8}{c}{$r=\frac{2-\alpha}{\alpha}$}\\[1.5pt]\cline{2-9}\rule{0pt}{9pt}

$\alpha=0.3\;\;$
&1.26e-3	&4.10e-4	&1.32e-4	&4.21e-5
    	&1.18e-2	&2.82e-3	&7.06e-4	&1.76e-4\\
&1.62	&1.64	&1.65	&&2.06	&2.00	&2.01\\[3pt]
$\alpha=0.5$
&3.26e-4	&1.03e-4	&3.32e-5	&1.10e-5
    	&7.87e-3	&2.01e-3	&4.98e-4	&1.24e-4\\
&1.67	&1.63	&1.59	&&1.97	&2.01	&2.00\\[3pt]	
$\alpha=0.7$
&6.77e-4	&2.58e-4	&1.01e-4	&4.02e-5
	   &5.64e-3	   &1.41e-3	   &3.48e-4	   &8.63e-5\\
&1.39	&1.35	&1.33   &&2.00	&2.02	&2.01\\
       	
\hline
      	
\end{tabular}}
\end{center}
\end{table}

\begin{table}[t!]
\begin{center}
\caption{%
Global maximum  nodal errors for $t\in[0,1]$
(odd rows) and
computational rates $q$ in $M^{-q}$ or $N^{-q}$ (even rows) on the graded mesh with $r=(2-\alpha)/\alpha$, $r=1$,
$r=2-\alpha$}
\label{t_global}
\vspace{-0.1cm}
\tabcolsep=4pt
{\small
\begin{tabular}{lrrrrrrrr}
\hline
\strut\rule{0pt}{10pt}
&\multicolumn{4}{l}{errors and convergence rates in time}&
\multicolumn{4}{l}{~~~errors and convergence rates in space}\\[0pt]
\strut\rule{0pt}{9pt}
&\multicolumn{4}{l}{$r=\frac{2-\alpha}{\alpha},\;\;\;N = \frac12M$}&
\multicolumn{4}{l}{~~~$r=\frac{2-\alpha}{\alpha},\;\;\;M = N^2
$}\\[2pt]
\hline
\strut\rule{0pt}{9pt}&
$M=2^{8}$& $M=2^{9}$& $M=2^{10}\!\!$& $M=2^{11}\!\!$&{}~~~$N=2^3$&~$N=2^4$&~$N=2^5$&~$N=2^6$\\
\hline
$\alpha=0.3\;\;$

&1.49e-4	&4.79e-5	&1.55e-5	&4.97e-6
       &1.96e-2	&4.82e-3	&1.20e-3	&3.01e-4\\
&1.64	&1.63	&1.64
        &&2.02	&2.00	&2.00	\\[3pt]

$\alpha=0.5$

&3.91e-4	&1.43e-4	&5.20e-5	&1.88e-5		
        &1.24e-2	&3.18e-3	&7.95e-4	&1.98e-4\\
&1.45	&1.46	&1.47	
        &&1.97	&2.00	&2.01\\[3pt]

$\alpha=0.7$

&8.90e-4	&3.83e-4	&1.63e-4	&6.83e-5
        &1.43e-2	&3.63e-3	&8.76e-4	&2.12e-4\\	
&1.22	&1.24	&1.25		
        &&1.98	&2.05	&2.05\\		
	        	
\hline
\strut\rule{0pt}{10pt}
&\multicolumn{4}{l}{errors and convergence rates in time}&
\\[0pt]
\strut\rule{0pt}{9pt}
&\multicolumn{4}{l}{$r=1,\;\;\;N = \frac1{128}M$}&
\multicolumn{4}{l}{~~~$r=2-\alpha,\;\;\;N = \frac14M$}\\[2pt]
\hline
\strut\rule{0pt}{9pt}&
$M=2^{15}\!\!\!$& $M=2^{16}\!\!\!$& $M=2^{17}\!\!\!$& $M=2^{18}\!\!\!$&{}~~~$M=2^{10}\!\!\!$& $M=2^{11}\!\!\!$& $M=2^{12}\!\!\!$& $M=2^{13}\!\!\!$\\
\hline
$\alpha=0.3\;\;$

&1.30e-2	&1.13e-2	&9.77e-3	&8.37e-3
       &9.77e-3	&7.47e-3	&5.59e-3	&4.11e-3\\
&0.20	&0.21	&0.22
        &&0.39	&0.42	&0.45\\[3pt]

$\alpha=0.5$

&2.73e-3	&1.95e-3	&1.39e-3	&9.88e-4		
        &2.73e-3	&1.64e-3	&9.88e-4	&5.94e-4\\
&0.49	&0.49	&0.49	
        &&0.73	&0.73	&0.73\\[3pt]

$\alpha=0.7$

&3.15e-4	&1.93e-4	&1.19e-4	&7.33e-5
        &9.84e-4	&5.27e-4	&2.82e-4	&1.51e-4\\	
&0.70	&0.70	&0.70	
        &&0.90	&0.90	&0.90\\		
	        	
\hline
\end{tabular}}
\end{center}
\end{table}

As a test problem, consider  \eqref{problem} with $\LL=-(\pt_{x_1}^2+\pt_{x_2}^2)$
and an Allen-Cahn type nonlinearity $f=(u^3-u)/\alpha$, posed
in  the square spatial domain $\Omega=(0,\pi)^2$ for $t\in[0,1]$,
subject to the initial condition
$u(0,t)=u_0=\frac25(2y-x^2)\,\sin x\,\sin y$.
We shall test the error bound
\eqref{FD_error_bound} of Theorem~\ref{theo_FD}(i) given for finite difference discretizations in space combined with the L1 scheme in time.
 The graded temporal mesh
$\{t_j=(j/M)^r\}_{j=0}^M$ will be used in all experiments.
 The spatial mesh is a uniform tensor product mesh of size $h=\pi/N$ (i.e. with $N$ equal mesh intervals in each coordinate direction).
 As the exact solution is unknown, the errors are computed using  the two-mesh principle.

 First, note that
 condition A2 is satisfied with $-\sigma_1=\sigma_2=1$, while
 the initial condition is in $[\sigma_1,\sigma_2]=[-1,1]$.
 In full agreement with Theorem~\ref{theo_FD}(ii), we have observed that  all our computed solutions were also in this range.

   Next, we look  into the more interesting case of convergence in positive time $t\gtrsim 1$ and give,
    in Table~\ref{t_positive_time}, the maximum nodal errors for
  the graded temporal meshes with  $r=1$, $r=(2-\alpha)/0.9$ and $r=(2-\alpha)/\alpha$.
  Recalling Remark~\ref{rem_positive_time},
  for $r=1$  we expect  convergence rates in time close to  $1$. The other two values
satisfy $r>2-\alpha$, for which our error bound~\eqref{FD_error_bound} combined with Remark~\ref{rem_positive_time} predicts the optimal convergence rate of order $2-\alpha$ with respect to time.
This clearly agrees with the computational convergence rates given in Table~\ref{t_positive_time}. The spatial convergence rates are close to $2$, which is also consistent with our theoretical bound.

The global maximum nodal errors for $t\in[0,1]$ were computed for the optimal grading parameter $r=(2-\alpha)/\alpha$ (see the upper part of Table~\ref{t_global}),  as well as for $r=1$ and $r=2-\alpha$ (see the lower part of the same table).
In view of Remark~\ref{rem_global_time},
 the theoretical error bound~\eqref{FD_error_bound}
 predicts
  the global convergence rates in time close to $\alpha r$, which is also in good agreement with the computational convergence rates in Table~\ref{t_global}.

Overall, we conclude that our numerical results are consistent with our theoretical findings.
We also refer the reader to numerical results in \cite{NK_XM}, which illustrate (for the linear case) that
our error bounds are remarkably sharp in the pointwise-in-time sense.

\end{document}